\documentclass[11pt]{article}
\usepackage{latexsym}

\newcommand{\eproof}{\mbox{\ }\hfill $\Box$ \par \vskip 10pt}

\newtheorem{Theorem}{Theorem}[section]
\newtheorem{lemma}[Theorem]{Lemma}
\newtheorem{prop}[Theorem]{Proposition}

\newtheorem{corol}[Theorem]{Corollary}

\begin{document}

\title{Transmission eigenvalue-free regions}

\author{{\sc Georgi Vodev}}

\date{}

\maketitle

\noindent
{\bf Abstract.} We prove the existence of large regions free of eigenvalues of the interior transmission problem.

\setcounter{section}{0}
\section{Introduction and statement of results}

Let $\Omega\subset{\bf R}^d$, $d\ge 2$, be a bounded, connected domain with a $C^\infty$ smooth boundary $\Gamma=\partial\Omega$. 
A complex number $\lambda\in {\bf C}$, $\lambda\neq 0$, will be said to be a transmission eigenvalue if the following problem has a non-trivial solution:
$$\left\{
\begin{array}{lll}
\left(\nabla c_1(x)\nabla+\lambda n_1(x)\right)u_1=0 &\mbox{in} &\Omega,\\
\left(\nabla c_2(x)\nabla+\lambda n_2(x)\right)u_2=0 &\mbox{in} &\Omega,\\
u_1=u_2,\,\,\, c_1\partial_\nu u_1=c_2\partial_\nu u_2& \mbox{on}& \Gamma,
\end{array}
\right.
\eqno{(1.1)}
$$
where $\nu$ denotes the exterior Euclidean unit normal to $\Gamma$, $c_j,n_j\in C^\infty(\overline\Omega)$, $j=1,2$ are strictly positive real-valued functions.
The purpose of this work is to study the localization of the possible transmission eigenvalues on ${\bf C}$ as $|\lambda|\to\infty$ under the condition
$$c_1(x)n_1(x)\neq c_2(x)n_2(x),\quad\forall x\in\Gamma.\eqno{(1.2)}$$
 Our first result is the following

\begin{Theorem}  Assume (1.2) together with the condition
$$c_1(x)=c_2(x),\quad \partial_\nu c_1(x)=\partial_\nu c_2(x),\quad\forall x\in\Gamma.\eqno{(1.3)}$$
Then there are no transmission eigenvalues in $\Lambda^-\cup\Lambda^+$, where
$$\Lambda^+:=\left\{\lambda\in{\bf C}:\,{\rm Re}\,\lambda\ge 0,\,|{\rm Im}\,\lambda|\ge 
C_\varepsilon\left({\rm Re}\,\lambda+1\right)^{\frac{3}{4}+\varepsilon}\right\}$$
for every $0<\varepsilon\ll 1$ with $C_\varepsilon>0$,
$$\Lambda^-:=\left\{\lambda\in{\bf C}:{\rm Re}\,\lambda\le -\widetilde C\right\}\cup
\left\{\lambda\in{\bf C}:-\widetilde C\le{\rm Re}\,\lambda\le 0,
\,\,|{\rm Im}\,\lambda|\ge C\right\}$$  $C,\widetilde C>0$ being  constants.
\end{Theorem}

\noindent
{\bf Remark 1.} In the case $c_1\equiv c_2\equiv 1$ it was proved previously in \cite{kn:LV3} that outside any sector $|{\rm Im}\,\lambda|\le \theta{\rm Re}\,\lambda$, $\forall\theta>0$,
there are at most a finite number of transmission eigenvalues. In the case $c_1\equiv c_2\equiv 1$, $n_1\equiv 1$ and $n_2(x)>1$, $\forall x\in\overline\Omega$,
the above theorem is proved in \cite{kn:H} but with $\Lambda^+$ replaced by a smaller eigenvalue-free region of the form
$$\left\{\lambda\in{\bf C}:\,{\rm Re}\,\lambda\ge 0,\,|{\rm Im}\,\lambda|\ge 
C\left({\rm Re}\,\lambda+1\right)^{\frac{24}{25}}\right\}.$$

The situation is far more interesting and different when the condition (1.3) is not fulfilled. In this case we have the following

\begin{Theorem} Assume (1.2) together with the condition
$$c_1(x)\neq c_2(x),\quad\forall x\in\Gamma.\eqno{(1.4)}$$
Then, there are no transmission eigenvalues in 
$$\left\{\lambda\in{\bf C}:{\rm Re}\,\lambda\ge 0,\,\,|{\rm Im}\,\lambda|\ge C\left({\rm Re}\,\lambda+1\right)^{\frac{4}{5}}\right\}$$
with a constant $C>0$. Moreover, if in addition we assume either the condition
$$\frac{n_1(x)}{c_1(x)}\neq \frac{n_2(x)}{c_2(x)},\quad\forall x\in\Gamma,\eqno{(1.5)}$$
or the condition
$$\frac{n_1(x)}{c_1(x)}=\frac{n_2(x)}{c_2(x)},\quad\forall x\in\Gamma,\eqno{(1.6)}$$
then there are no transmission eigenvalues in $\Lambda^+$. Under the condition
$$(c_1(x)-c_2(x))(c_1(x)n_1(x)-c_2(x)n_2(x))> 0,\quad\forall x\in\Gamma,\eqno{(1.7)}$$
 there are no transmission eigenvalues in $\Lambda^-$. Finally, if we assume the condition
$$(c_1(x)-c_2(x))(c_1(x)n_1(x)-c_2(x)n_2(x))< 0,\quad\forall x\in\Gamma,\eqno{(1.8)}$$
then for every $N\ge 1$ there is a constant $C_N>0$ so that there are no transmission eigenvalues in 
 $$\left\{\lambda\in{\bf C}:{\rm Re}\,\lambda\le 0,\,\,|{\rm Im}\,\lambda|\ge C_N\left(|{\rm Re}\,\lambda|+1\right)^{-N}\right\}.$$ 
\end{Theorem}

\noindent
{\bf Remark 2.} Note that the condition (1.8) implies (1.5).

\noindent
{\bf Remark 3.} It is clear from the proof that the fact that we can take an arbitray $N$ above comes from the $C^\infty$- smoothness of
the boundary $\Gamma$ and the coefficients $c_j$, $n_j$ near $\Gamma$. Therefore, it is natural to expect that if more regularity is assumed
(e.g. Gevrey class or analyticity), a larger eigenvalue-free region exists. Indeed, using the techniques of \cite{kn:Sj}
one can show that in the analytic case there is a region free of eigenvalues of the form
$$\left\{\lambda\in{\bf C}:{\rm Re}\,\lambda\le 0,\,\,|{\rm Im}\,\lambda|\ge C\exp\left(-\beta|{\rm Re}\,\lambda|^{1/2}\right)\right\}$$
with some constants $C,\beta>0$. 

\noindent
{\bf Remark 4.} It is clear from our construction of the parametrix in the region ${\rm Re}\,\lambda<0$ that under the condition (1.8)
one can construct quasimodes for the problem (1.1) concentrated in an arbitrary neighbourhood of the boundary $\Gamma$ due to the existence
in this case of surface waves moving on $\Gamma$ with a speed $\sqrt{c(x)}$, where $c$ denotes the restriction on $\Gamma$ of the
function $|c_1^2-c_2^2|/|c_1n_1-c_2n_2|$. These waves are very similar to the Rayleigh surface waves in the linear elasticity
studied in \cite{kn:SjV}, \cite{kn:SV1}, \cite{kn:SV2} and have practically the same properties. In particular, as in \cite{kn:SV1}
one can show that these quasimodes imply the existence of infinitely many transmission eigenvalues with negative real parts.
In fact, much more can be proved, namely an asymptotic of the counting function $N^-(r)=\#\{\lambda-{\rm trans.\, eig.}: 
{\rm Re}\,\lambda<0,\,|\lambda|\le r^2\}$, $r>1$. Indeed, as in \cite{kn:SjV} one can show that
$$N^-(r)=\left(\frac{r}{2\pi}\right)^{d-1}\omega_{d-1}\int_\Gamma c(x)^{-\frac{d-1}{2}}dx+O(r^{d-2})$$
where $\omega_{d-1}:={\rm Vol}\{x\in{\bf R}^{d-1}:|x|\le 1\}$, 
provided the multiplicity of an eigenvalue $\lambda_k$ is defined by
$${\rm mult}(\lambda_k)={\rm tr}\,(2i\pi)^{-1}\int_{|\lambda-\lambda_k|=\varepsilon}\left(c_1\frac{d{\cal N}_1}{d\lambda}(\lambda)-
c_2\frac{d{\cal N}_2}{d\lambda}(\lambda)\right)\left(c_1{\cal N}_1(\lambda)-c_2{\cal N}_2(\lambda)\right)^{-1}d\lambda,$$
 $0<\varepsilon\ll 1$, where ${\cal N}_j(\lambda)$ denotes the Dirichlet-to-Neumann map corresponding to the pair $(n_j,c_j)$.
 
 \noindent
{\bf Remark 5.} These kind of eigenvalue-free regions are crucial for bounding the remainder in the asymptotics of the counting
function of all transmission eigenvalues. Indeed, it has been proved recently in \cite{kn:PV} that the total counting function
$N(r)=\#\{\lambda-{\rm trans.\, eig.}:\,|\lambda|\le r^2\}$, $r>1$, satisfies the asymptotics
$$N(r)=(\tau_1+\tau_2)r^d+O_\varepsilon(r^{d-\kappa+\varepsilon}),\quad\forall\,0<\varepsilon\ll 1,$$
where $0<\kappa\le 1$ is such that there are no transmission eigenvalues in the region
$$\left\{\lambda\in{\bf C}:\,|{\rm Im}\,\lambda|\ge C\left(|{\rm Re}\,\lambda|+1\right)^{1-\frac{\kappa}{2}}\right\},\quad C>0,$$
and $$\tau_j=\frac{\omega_d}{(2\pi)^d}\int_\Omega\left(\frac{n_j(x)}{c_j(x)}\right)^{d/2}dx,$$
$\omega_d$ being the volume of the unit ball in ${\bf R}^d$.

\begin{corol}  Assume (1.2) together with the condition
$$n_1(x)=n_2(x),\quad\forall x\in\Gamma.\eqno{(1.9)}$$
Then there are no transmission eigenvalues in $\Lambda^-\cup\Lambda^+$.
\end{corol}

There have been recently many works studing mainly the discreteness and the asymptotic behaviour of the counting function of the transmission
eigenvalues (see \cite{kn:DP}, \cite{kn:H}, \cite{kn:LV1}, \cite{kn:LV2}, \cite{kn:LV3}, \cite{kn:LV4}, \cite{kn:PS}, \cite{kn:R1}, \cite{kn:R2},
 \cite{kn:S} and the references therein). For example, in \cite{kn:LV1} Weyl type asymptotics have been proved for the counting
function of all transmission eigenvalues under the condition (1.4), while in \cite{kn:LV4} a lower bound of right order has been proved
for the counting function of the transmission eigenvalues belonging to $(0,+\infty)$. Under the conditions (1.2) and (1.3), upper bounds of 
the counting function of all transmission eigenvalues have been proved in \cite{kn:DP}, \cite{kn:R1}, \cite{kn:R2}, and an asymptotic has been
given in \cite{kn:PS} when $\Omega$ is a ball and the coefficients are constants. 

To prove the eigenvalue-free regions we transform the problem (1.1) into a semi-classical one by putting $h=|{\rm Re}\,\lambda|^{-1/2}$,
$z=\frac{\lambda}{|{\rm Re}\,\lambda|}$, if $|{\rm Re}\,\lambda|\ge |{\rm Im}\,\lambda|$, $|{\rm Re}\,\lambda|\gg 1$, and 
$h=|{\rm Im}\,\lambda|^{-1/2}$,
$z=\frac{\lambda}{|{\rm Im}\,\lambda|}$, if $|{\rm Im}\,\lambda|\ge |{\rm Re}\,\lambda|$, $|{\rm Im}\,\lambda|\gg 1$.
Thus we have to show that if $h$ is small enough and $z/h^2$ belongs to the eigenvalue-free regions described above, then under the corresponding conditions the solutions
to (1.1) are identically zero. In fact, it suffices to show that $u_1|_\Gamma$ is identically zero since this would imply that $u_1$ and $u_2$
are identically zero, too. To do so, we construct in Section 3 a parametrix of the solutions to the interior boundary value problem (see equation (3.1) below)
near the boundary $\Gamma$ by using $h$-FIOs with a complex-valued phase satisfying the eikonal equation mod $O(x_1^N)$, $N>1$ being an arbitrary
integer and 
$0<x_1\ll 1$ is the normal coordinate to the boundary (which is nothing else but the Euclidean distantce from a point $x\in\Omega$ to $\Gamma$). 
The amplitude must satisfy the transport equations mod $O(x_1^N)$. We solve these equations in Section 4. Furthermore, we use this parametrix
to show that the Dirichlet-to-Neumann map can be approximated by $h-\Psi $DOs belonging (uniformly in $z$) to the class ${\cal S}^1_0$ (sse Section 2 for the definition) if ${\rm Re}\,z=-1$, $|{\rm Im}\,z|\le 1$ or $|{\rm Re}\,z|\le 1$, $|{\rm Im}\,z|=1$, and to 
${\cal S}^1_{1/2-\epsilon}$ if ${\rm Re}\,z=1$, $h^{\frac{1}{2}-\epsilon}\le |{\rm Im}\,z|\le 1$, $0<\epsilon\ll 1$. Thus we reduce the problem
of finding eigenvalue-free regions to that one of inverting $h-\Psi $DOs (depending on an additional parameter $z$) on a compact manifold. Note that these classes of $h-\Psi $DOs
are nice in the sense that there is a symbol calculas for them as well as a simple criteria of $L^2$ boundeness (e.g. see \cite{kn:DS}).
We recall these properties of the $h-\Psi $DOs in Section 2. 
In particular, to invert such an operator it suffices to invert its principal symbol and determine the class of symbols the inverse belongs to.
That is precisely what we do in Section 5. Note that the study of the case $n=1$ in \cite{kn:S} (see also \cite{kn:PS}) suggests that 
there are probably larger eigenvalue-free regions in ${\rm Re}\,\lambda>0$ at least in some specific cases as for example $\Omega$ is
a ball and $c_j,n_j$ constants. In this latter case one has to invert Bessel functions instead of $h-\Psi $DOs, which seems to be much easier.
In the general case studied here, however, it would be impossible to do better since it is impossible to construct a parametrix for the equation (3.1)
when ${\rm Re}\,z=1$, $0<|{\rm Im}\,z|\ll h^{\frac{1}{2}-\epsilon}$. As a consequence, in this region the Dirichlet-to-Neumann map
is no longer an $h-\Psi $DO, and hence it is impossible to use the theory of the $h-\Psi $DOs to invert our operator. 

\section{$h$-pseudo-differential operators on a compact manifold}

Let $X$ be a $C^\infty$ smooth compact manifold without boundary, $n={\rm dim}\,X\ge 1$.
Let $(x,\xi)$ be coordinates on $T^*X$ and let $a\in C^\infty(T^*X)$. Then the $h$-pseudo-differential operator with a symbol $a$ is defined as follows
$$\left({\rm Op}_h(a)f\right)(x)=\left(\frac{1}{2\pi h}\right)^n\int_{T^*X}e^{-\frac{i}{h}\langle x-y,\xi\rangle}a(x,\xi)f(y)dyd\xi$$
where $h>0$ is a small parameter. Of course, in order that this operator has {\it nice} properties the function $a$ must belong to some
class of symbols. In what follows in this section we will introduce several classes of symbols which will play important role in our analysis. 
First, given $\ell\in{\bf R}$, $\delta_1,\delta_2\ge 0$ and a function $\mu>0$, we denote
by $S^{\ell}_{\delta_1,\delta_2}(\mu)$ the set of all functions $a\in C^\infty(T^*X)$ such that 
$$\left|\partial_x^\alpha\partial_\xi^\beta a(x,\xi)\right|\le C_{\alpha,\beta}\mu^{\ell-\delta_1|\alpha|-\delta_2|\beta|}$$
for all multi-indices $\alpha,\beta$ with constants $C_{\alpha,\beta}>0$ independent of $h$, $\mu$.
The following simple properties will be often used in the next sections: If $a\in S^{\ell}_{\delta_1,\delta_2}(\mu)$,
then $\partial_x^\alpha\partial_\xi^\beta a\in S^{\ell-\delta_1|\alpha|-\delta_2|\beta|}_{\delta_1,\delta_2}(\mu)$.
If $a_j\in S^{\ell_j}_{\delta_1,\delta_2}(\mu)$, $j=1,2$,
then $a_1a_2\in S^{\ell_1+\ell_2}_{\delta_1,\delta_2}(\mu)$. If $b(x)\in C^\infty(X)$, independent of $\xi$, and 
$a\in S^{\ell}_{\delta_1,\delta_2}(\mu)$, then $ba\in S^{\ell}_{\delta_1,\delta_2}(\mu)$ if $\mu\le Const$ or $\mu\ge Const>0$ and $\delta_1=0$.
We also need a simple criteria for this class of operators to be bounded on $L^2(X)$.

\begin{prop} Let the function $a$ satisfy
$$\sup_{x,\xi\in T^*X}\left|\partial_x^\alpha a(x,\xi)\right|=C_{\alpha}<\infty\eqno{(2.1)}$$
for all multi-indices $\alpha$. Then the operator ${\rm Op}_h(a)$ is bounded on $L^2(X)$
and
$$\left\|{\rm Op}_h(a)\right\|_{L^2(X)\to L^2(X)}\le C\sum_{|\alpha|\le n+1}C_\alpha h^{|\alpha|/2}\eqno{(2.2)}$$
with a constant $C>0$ independent of $h$ and $C_\alpha$.
In particular, if $a\in S^{\ell}_{\delta,\delta_2}(\mu)$ with $\ell\le 0$ and $\mu(x,\xi)\ge\mu_0>0$, we have the bound 
$$\left\|{\rm Op}_h(a)\right\|_{L^2(X)\to L^2(X)}\le C\mu_0^{\ell}\left(1+ \frac{\sqrt{h}}{\mu_0^{\delta}}\right)^{n+1}\eqno{(2.3)}$$
with a constant $C>0$ independent of $h$ and $\mu_0$. 
 \end{prop}

{\it Proof.} It is based on the observation that the boundness of ${\rm Op}_h(a)$ on $L^2$ is equivalent to that of
the classical operator ${\rm Op}_1(a_h)$, where $a_h(x,\xi)=a(\sqrt{h}x,\sqrt{h}\xi)$. On the other hand, since $X$ is compact, it is well known
(see Theorem 18.1.11$'$ of \cite{kn:Ho}) 
that the norm of ${\rm Op}_1(a_h):L^2\to L^2$ is bounded by  $\sum_{|\alpha|\le n+1}\sup|\partial_x^\alpha a_h(x,\xi)|$,
which implies (2.2).
\eproof

Given $k\in{\bf R}$, $0\le\delta\le\frac{1}{2}$, denote by ${\cal S}^{k}_\delta$ the set of all functions $a\in C^\infty(T^*X)$ satisfying
$$\left|\partial_x^\alpha\partial_\xi^\beta a(x,\xi)\right|\le C_{\alpha,\beta}h^{-\delta(|\alpha|+|\beta|)}
\langle\xi\rangle^{k-|\beta|}$$
for all multi-indices $\alpha,\beta$ with constants $C_{\alpha,\beta}>0$ independent of $h$. We will denote by ${\rm OP}{\cal S}^{k}_\delta$
the set of the $h$-pseudo-differential operators with symbols in ${\cal S}^{k}_\delta$. 
It follows from the above proposition
that if $a\in {\cal S}^0_\delta$, then ${\rm Op}_h(a):L^2\to L^2=O(1)$. It is also well-known (e.g. see Section 7 of \cite{kn:DS}) that when $\delta<\frac{1}{2}$
there is a nice symbol calculas and in particular the symbol of the composition of $h$-pseudo-differential operators with symbols in this class
can be calculated explicitly mod $O(h^\infty)$. Thus, if $a\in{\cal S}^{k}_\delta$ with $0\le\delta<\frac{1}{2}$ and 
$|a|\ge C\langle\xi\rangle^{k}$ with $C>0$ independent of $h$, then the operator ${\rm Op}_h(a)$ is invertible with an inverse
belonging to ${\rm OP}{\cal S}^{-k}_\delta$.
 The following proposition is essentially proved in Section 7 of \cite{kn:DS}.
Here we sketch the proof for the sake of completeness.

\begin{prop} Let $h^{\ell_\pm}a^\pm\in{\cal S}^{\pm k}_{\delta}$, $\delta<\frac{1}{2}$, where ${\ell_\pm}\ge 0$ are some numbers. 
Assume in addition that the functions $a^\pm$ satisfy
$$\left|\partial_x^{\alpha_1}\partial_\xi^{\beta_1} a^+(x,\xi)\partial_x^{\alpha_2}\partial_\xi^{\beta_2}a^-(x,\xi)\right|\le \mu_0C_{\alpha_1,\beta_1,\alpha_2,\beta_2}h^{-\frac{|\alpha_1|+|\beta_1|+|\alpha_2|+|\beta_2|}{2}}\eqno{(2.4)}$$
for all multi-indices $\alpha_1,\beta_1,\alpha_2,\beta_2$ such that $|\alpha_j|+|\beta_j|\ge 1$, $j=1,2$, with constants $C_{\alpha_1,\beta_1,\alpha_2,\beta_2}>0$ independent of $h$, $\mu_0$.
Then we have
 $$\left\|{\rm Op}_h(a^+){\rm Op}_h(a^-)-{\rm Op}_h(a^+a^-)\right\|_{L^2(X)\to L^2(X)}\le 
 C\mu_0+Ch\eqno{(2.5)}$$
 with a constant $C>0$ independent of $h$ and $\mu_0$.
\end{prop}

{\it Proof.} In view of formula (7.15) 
of \cite{kn:DS} the operator in the left-hand side of (2.5)
whose norm we would like to bound is an $h$-psdo with symbol $b(x,\xi,x,\xi)$, where the function $b$ is given by
$$b(x,\xi,y,\eta)=\left(e^{ih D_\xi\cdot D_y}-1\right)a(x,\xi,y,\eta)$$
where we have put $a=a^+(x,\xi)a^-(y,\eta)$ and $D=-i\partial$. 
It follows from the analysis in Section 7 of \cite{kn:DS} (see (7.17) and (7.19)) that given any integer $N\ge 2$ the function $b$ can be decomposed as $b_N+\widetilde b_N$, where
$$b_N=\sum_{j=1}^{N-1}\frac{1}{j!}(ih D_\xi\cdot D_y)^ja
=\sum_{j=1}^{N-1}\frac{(ih)^j}{j!}\sum_{|\alpha|=j}D_\xi^\alpha a^+(x,\xi)D_y^\alpha a^-(y,\eta)$$
while the remainder $\widetilde b_N$ satisfies
$$\left|\partial_x^\alpha\partial_y^\beta\widetilde b_N(x,\xi,y,\eta)\right|\le C_{\alpha,\beta}h^{N(1-2\delta)-\ell}\langle\xi\rangle^{k}\langle\eta\rangle^{-k}
\le C_{\alpha,\beta}h^{N(1-2\delta)-\ell}$$
if $\eta=\xi$, where $\ell=\ell_++\ell_-+s_n+\delta(|\alpha|+|\beta|)$ is independent of $N$. In view of Proposition 2.1, this implies
that there exists some $\ell_1>0$ independent of $N$ such that
$$\left\|{\rm Op}_h(\widetilde b_N(x,\xi,x,\xi))\right\|_{L^2\to L^2}\le Ch^{N(1-2\delta)-\ell_1}\le Ch\eqno{(2.6)}$$
if $N$ is taken large enough. On the other hand, it is easy to see that (2.4) implies
$$\left|\partial_x^\alpha\partial_y^\beta b_N(x,\xi,x,\xi)\right|\le \mu_0C_{\alpha,\beta}h^{-\frac{|\alpha|+|\beta|}{2}}.$$
 By Proposition 2.1,
$$\left\|{\rm Op}_h(b_N(x,\xi,x,\xi))\right\|_{L^2\to L^2}\le C\mu_0.\eqno{(2.7)}$$
Clearly, (2.5) follows from (2.6) and (2.7).
\eproof

\section{Parametrix near the boundary}

Let $z\in Z=Z_1\cup Z_2\cup Z_3$, where $Z_1=\{z\in{\bf C}:{\rm Re}\,z=1,\,0<|{\rm Im}\,z|\le 1\}$,
$Z_2=\{z\in{\bf C}:{\rm Re}\,z=-1,\,|{\rm Im}\,z|\le 1\}$,
$Z_3=\{z\in{\bf C}:|{\rm Re}\,z|\le 1,\,|{\rm Im}\,z|=1\}$.  Clearly, we have $1\le|z|\le 2$. 
Given any $f\in L^2(\Gamma)$ let $u$ solve the equation
$$\left\{
\begin{array}{lll}
 \left(P(h)-z\right)u=0&\mbox{in}&\Omega,\\
 u=f&\mbox{on}&\Gamma,
\end{array}
\right.
\eqno{(3.1)}
$$
where 
$$P(h)=-\frac{h^2}{n(x)}\nabla c(x)\nabla$$
and $h>0$ is a small parameter, $c,n\in C^\infty(\overline\Omega)$ being strictly positive functions.
Let $(x',\xi')$ be coordinates on $T^*\Gamma$ and denote by $r_0(x',\xi')$ the principal symbol of the Laplace-Beltrami operator,
$-\Delta_\Gamma$, on $\Gamma$ equipped with the Riemannian metric induced by the Euclidean one in ${\bf R}^d$. It is well-known that
$r_0$ is a polynomial function in $\xi'$, homogeneous of order 2, and $C_2|\xi'|^2\ge r_0(x',\xi')\ge C_1|\xi'|^2$ with constants $C_2>C_1>0$.
Set $m(x)=\frac{n(x)}{c(x)}$ and denote by $\gamma$ the restriction on $\Gamma$, that is, $\gamma m=m|_{\Gamma}$. 
Define the function $\rho\in C^\infty(T^*\Gamma)$ as being the root of the equation
$$\rho^2+r_0(x',\xi')-\gamma m(x')z=0$$
with ${\rm Im}\,\rho>0$ (which is easily seen to exist as long as $z\in Z$). 
  In what follows in this paper $C$ and $\widetilde C$ will denote  positive
constants independent of $z$, $h$ and $f$, which may change from line to line.

\begin{lemma} Let $z\in Z_1\cup Z_3$. Then 
$${\rm Im}\,\rho\ge \frac{|{\rm Im}\,z|}{2|\rho|},\eqno{(3.2)}$$
$$|\rho|\ge C\sqrt{|{\rm Im}\,z|},\eqno{(3.3)}$$
while for $r_0\ge 2\gamma m$, we have
$$\widetilde C\sqrt{r_0+1}\ge 2{\rm Im}\,\rho\ge|\rho|\ge C\sqrt{r_0+1}.\eqno{(3.4)}$$
 Let $z\in Z_2$. Then (3.4) holds for all $r_0\ge 0$.
\end{lemma}

{\it Proof.} Clearly, (3.2) follows from the identity
$$2{\rm Im}\,\rho{\rm Re}\,\rho=\gamma m{\rm Im}\,z.$$
The bound (3.3) follows easily from the identity
$$|\rho|^4=|r_0-\gamma mz|^2=(r_0-\gamma m{\rm Re}\,z)^2+(\gamma m{\rm Im}\,z)^2.$$
For $r_0\ge 2\gamma m$, we have
$$|\rho|^4\ge \frac{1}{4}(r_0+\gamma m{\rm Re}\,z)^2+(\gamma m{\rm Im}\,z)^2$$
and
$$({\rm Im}\,\rho)^2-({\rm Re}\,\rho)^2=r_0-\gamma m{\rm Re}\,z\ge 0$$
 When $z\in Z_2$ these inequalities clearly hold for all $r_0\ge 0$.
\eproof

Let $\phi\in C^\infty({\bf R})$,
$\phi(\sigma)=1$ for $|\sigma|\le 1$, $\phi(\sigma)=0$ for $|\sigma|\ge 2$, and set
$$\chi(x',\xi')=\phi\left(\delta_0r_0(x',\xi')\right)$$
where $0<2\delta_0\le \min_{x'\in\Gamma}\frac{1}{\gamma m(x')}$. 
We will say that a function $a\in C^\infty(T^*\Gamma)$ belongs to $ S^{\ell_1}_{\delta_1,\delta_2}(\mu_1)+S^{\ell_2}_{\delta_3,\delta_4}(\mu_2)$
if $\chi a\in S^{\ell_1}_{\delta_1,\delta_2}(\mu_1)$ and $(1-\chi)a\in S^{\ell_2}_{\delta_3,\delta_4}(\mu_2)$.

\begin{lemma} We have $\rho,|\rho|\in S^{1}_{2,2}(|\rho|)+S^1_{0,1}(|\rho|)$,
$\rho^{-1},|\rho|^{-1}\in S^{-1}_{2,2}(|\rho|)+S^{-1}_{0,1}(|\rho|)$ uniformly in $z$.
\end{lemma}

{\it Proof.} In view of Lemma 3.1, we have $|\rho|\le C$ on supp$\,\chi$, $|\rho|\ge C$ on supp$\,(1-\chi)$, $C>0$. We have to show that the function $\rho$ satisfies the estimates
$$\left|\partial_{x'}^\alpha\partial_{\xi'}^\beta \rho\right|\le C_{\alpha,\beta}|\rho|^{1-2|\alpha|-2|\beta|}\quad
\mbox{on}\quad {\rm supp}\,\chi,\eqno{(3.5)}$$
$$\left|\partial_{x'}^\alpha\partial_{\xi'}^\beta \rho\right|\le C_{\alpha,\beta}|\rho|^{1-|\beta|}\quad
\mbox{on}\quad {\rm supp}\,(1-\chi),\eqno{(3.6)}$$
for all multi-indices $\alpha$ and $\beta$ with constants $C_{\alpha,\beta}>0$ independent of $z$, and similarly for the function $|\rho|$.
We will proceed by induction in $K=|\alpha|+|\beta|$. 
Differentiating the above equation we get
$$E_{\alpha,\beta}:=-\partial_{x'}^\alpha\partial_{\xi'}^\beta\left(r_0(x',\xi')-\gamma m(x')z\right) =
\partial_{x'}^\alpha\partial_{\xi'}^\beta(\rho^2)=2\rho\partial_{x'}^\alpha\partial_{\xi'}^\beta\rho+F_{\alpha,\beta}$$
where $F_{\alpha,\beta}$ is a linear combination of functions of the form $\partial_{x'}^{\alpha_1}\partial_{\xi'}^{\beta_1}\rho
\partial_{x'}^{\alpha_2}\partial_{\xi'}^{\beta_2}\rho$ with multi-indices satisfying $|\alpha_1|+|\alpha_2|=|\alpha|$, $|\beta_1|+|\beta_2|=
|\beta|$, $|\alpha_1|+|\beta_1|\le K-1$, $|\alpha_2|+|\beta_2|\le K-1$. Hence, assuming (3.5) and (3.6) fulfilled for 
$|\alpha|+|\beta|\le K-1$ leads to the consclusion that $F_{\alpha,\beta}=O\left(|\rho|^{2-2|\alpha|-2|\beta|}\right)$ on supp$\,\chi$
and $F_{\alpha,\beta}=O\left(|\rho|^{2-|\beta|}\right)$ on supp$\,(1-\chi)$. On the other hand, we have $E_{\alpha,\beta}=O(1)$ on supp$\,\chi$
and $E_{\alpha,\beta}=O\left(|\rho|^{2-|\beta|}\right)$ on supp$\,(1-\chi)$ uniformly in $z$. From this and the above identity we conclude that
(3.5) and (3.6) hold for $|\alpha|+|\beta|= K$, as desired. The proof concerning the function $|\rho|$ is similar, using the identity
$$\widetilde E_{\alpha,\beta}:=\partial_{x'}^\alpha\partial_{\xi'}^\beta\left((r_0(x',\xi')-\gamma m(x'){\rm Re}\,z)^2+(\gamma m(x'){\rm Im}\,z)^2\right)$$ $$ =
\partial_{x'}^\alpha\partial_{\xi'}^\beta(|\rho|^4)=4|\rho|^3\partial_{x'}^\alpha\partial_{\xi'}^\beta|\rho|+\widetilde F_{\alpha,\beta}$$
where $\widetilde F_{\alpha,\beta}$ is a linear combination of functions of the form $\partial_{x'}^{\alpha_1}\partial_{\xi'}^{\beta_1}|\rho|
\partial_{x'}^{\alpha_2}\partial_{\xi'}^{\beta_2}|\rho|\partial_{x'}^{\alpha_3}\partial_{\xi'}^{\beta_3}|\rho|
\partial_{x'}^{\alpha_4}\partial_{\xi'}^{\beta_4}|\rho|$ with multi-indices satisfying $|\alpha_1|+|\alpha_2|+|\alpha_3|+|\alpha_4|=|\alpha|$, $|\beta_1|+|\beta_2|+|\beta_3|+|\beta_4|=|\beta|$, $|\alpha_j|+|\beta_j|\le |\alpha|+|\beta|-1$, $j=1,2,3,4$.
Clearly, on supp$\,\chi$ we have $\widetilde E_{\alpha,\beta}=O(|\rho|^2)$ for $|\alpha|+|\beta|=1$, $\widetilde E_{\alpha,\beta}=O(1)$ for $|\alpha|+|\beta|\ge 2$, while on supp$\,(1-\chi)$ we have
 $\widetilde E_{\alpha,\beta}=O\left(|\rho|^{4-|\beta|}\right)$. Therefore, the estimates (3.5) and (3.6) for the function $|\rho|$ can be
 proved by induction in $|\alpha|+|\beta|$ as above. The function $\rho^{-1}$ (resp. $|\rho|^{-1}$) can be treated similarly using the identity
$$0=\partial_{x'}^\alpha\partial_{\xi'}^\beta\left(\rho\rho^{-1}\right)=\rho \partial_{x'}^\alpha\partial_{\xi'}^\beta\left(\rho^{-1}\right)
+{\cal F}_{\alpha,\beta}$$
for $|\alpha|+|\beta|\ge 1$, where ${\cal F}_{\alpha,\beta}$ is a linear combination of functions of the form $\partial_{x'}^{\alpha_1}\partial_{\xi'}^{\beta_1}\rho
\partial_{x'}^{\alpha_2}\partial_{\xi'}^{\beta_2}(\rho^{-1})$ with multi-indices satisfying $|\alpha_1|+|\alpha_2|=|\alpha|$, $|\beta_1|+|\beta_2|=
|\beta|$, $|\alpha_j|+|\beta_j|\le |\alpha|+|\beta|-1$, $j=1,2$.
\eproof

Denote ${\cal D}_\nu=-ih\partial_\nu$, $Z_{1,\varepsilon}:=\{z\in Z_1:|{\rm Im}\,z|\ge h^{\frac{1}{2}-\varepsilon}\}$, where
$0\le\varepsilon\ll 1$. We also equipe the Sobolev space $H^1(\Gamma)$ with the semi-classical norm
$\|f\|_{H^1(\Gamma)}=\sum_{|\alpha|\le 1}h^{|\alpha|}\|\partial_{x'}^\alpha f\|_{L^2(\Gamma)}$. 

\begin{Theorem} Given any $0<\epsilon\ll 1$ there is $0<h_0(\epsilon)\ll 1$ so that for $z\in Z_{1,\epsilon}$
and $0<h\le h_0$ the solution $u$ to (3.1) satisfies the estimate
$$\left\|\gamma{\cal D}_\nu u-{\rm Op}_h(\rho+hb)f\right\|_{H^1(\Gamma)}\le \frac{Ch}{\sqrt{|{\rm Im}\,z|}}
\|f\|_{L^2(\Gamma)}\eqno{(3.7)}$$
where $b\in S^{0}_{0,1}(\langle\xi'\rangle)$ does not depend on $h$, $z$ and the function $n$.
Moreover, (3.7) holds for all $z\in Z_2\cup Z_3$ with $|{\rm Im}\,z|$ replaced by 1.
\end{Theorem}

{\it Proof.} To prove (3.7) we will construct a parametrix to the solution of (3.1) near the boundary $\Gamma$. In fact, it suffices to
carry out this construction locally and then to glue up all pices by using a partition of the unity on $\Gamma$. Indeed, it is well-known that 
given an arbitrary point $x^0\in\Gamma$, there exists a small neighbourhood ${\cal O}(x^0)\subset\overline\Omega$ of $x^0$
and local coordinates $(x_1,x')\in {\cal O}(x^0)$ such that $x^0=(0,0)$, $\Gamma\cap {\cal O}(x^0)$ is defined by $x_1=0$,
$x'$ being coordinates in $\Gamma\cap {\cal O}(x^0)$, $x_1>0$ in $\Omega\cap {\cal O}(x^0)$, and in these coordinates the operator
$${\cal P}(z,h)=-\frac{h^2}{c(x)}\nabla c(x)\nabla-z\frac{n(x)}{c(x)}$$
can be written in the form
$${\cal P}(z,h)={\cal D}_{x_1}^2+r(x,{\cal D}_{x'})-zm(x) + hq(x,{\cal D}_x)+h^2\widetilde q(x),$$
where we have put ${\cal D}_{x_1}=-ih\partial_{x_1}$, ${\cal D}_{x'}=-ih\partial_{x'}$, $r(x,\xi')=\langle R(x)\xi',\xi'\rangle$, 
$R=(R_{ij})$ being a symmetric $(d-1)\times(d-1)$ matrix-valued function with smooth real-valued entries, $q(x,\xi)=\langle q(x),\xi\rangle$, 
$q(x)$ and $\widetilde q(x)$ being smooth functions.
Moreover, we have $r(0,x',\xi')=r_0(x',\xi')$, the principal symbol of $-\Delta_\Gamma$ written in the coordinates $(x',\xi')$.
Let $\psi(x')\in C_0^\infty(\Gamma\cap {\cal O}(x^0))$, $\psi=1$ in a neighbourhood of $x^0$. We will construct a parametrix, $\widetilde u_\psi$, 
of (3.1), $\widetilde u_\psi|_{x_1=0}=\psi f$, in the form
$$\widetilde u_\psi(x)=(2\pi h)^{-d+1}\int\int e^{\frac{i}{h}\varphi(x,y',\xi',z)}
\Phi_\delta(x,\xi',z)a(x,\xi',z,h)f(y')dy'd\xi',$$
where $\Phi_\delta=\phi\left(\frac{x_1}{\delta}\right)\phi\left(\frac{x_1}{\delta\rho_1}\right)$, $\rho_1=1$ if $z\in Z_2\cup Z_3$, $\rho_1=|\rho|^3$ if $z\in Z_1$, 
$\phi$ being as above, $\delta>0$ is a small constant independent of $x,\xi',h,z$ to be fixed later on. The phase $\varphi$ is a
complex-valued function such that $\varphi|_{x_1=0}=-\langle x'-y',\xi'\rangle$, and the amplitude $a$ satisfies $a|_{x_1=0}=\psi(x')$.
More generally, given any integer $N\gg 1$ we will be searching $\varphi$ and $a$ in the form
$$\varphi=-\langle x'-y',\xi'\rangle+\sum_{k=1}^{N-1} x_1^k\varphi_k(x',\xi',z),$$
$$a=\sum_{k=0}^{N-1} \sum_{j=0}^{N-1} x_1^kh^ja_{k,j}(x',\xi',z)$$
so that $\varphi$ satisfies the eikonal equation mod $O(x_1^N)$:
$$\left(\partial_{x_1}\varphi\right)^2+r(x,\nabla_{x'}\varphi)-m(x)z=x_1^N\Psi_N(x,\xi',z)\eqno{(3.8)}$$
 and $a$ satisfies the equation
$$e^{-\frac{i}{h}\varphi}{\cal P}(z,h)e^{\frac{i}{h}\varphi}a=x_1^NA_N(x,\xi',z,h)+h^NB_N(x,\xi',z,h)\eqno{(3.9)}$$
where $\Psi_N$, $A_N$ and $B_N$ are smooth functions.
In Section 4 we will prove the following

\begin{prop} Let $z\in Z_{1,0}\cup Z_2\cup Z_3$. Then, for a suitable choice of the constant $\delta$, the equations (3.8) and (3.9) have smooth
solutions $\varphi$ and $a$ of the form above, $\varphi=-\langle x'-y',\xi'\rangle+\widetilde\varphi$, with $\varphi_1=\rho$, 
$a_{0,0}=\psi$, $a_{0,j}=0$ for $j\ge 1$, $a_{1,j}\in S^{-1-4j}_{2,2}(|\rho|)+S^{-j}_{0,1}(|\rho|)$, $j\ge 0$,
$$a_{1,0}=-\frac{i}{2} q(0,x',1,\xi'/\rho)\psi-\frac{1}{2\rho}\langle R(0,x')\xi',\nabla_{x'}\psi(x')\rangle.$$
Moreover, for all integers $k\ge 0$, we have
$$x_1^{-1}\widetilde\varphi\in S^{1}_{2,2}(|\rho|)+S^{1}_{0,1}(|\rho|),\quad   
\partial_{x_1}^k\widetilde\varphi\in S^{4-3k}_{2,2}(|\rho|)+S^{1}_{0,1}(|\rho|),$$  
 $$\partial_{x_1}^ka\in S^{2-3k}_{2,2}(|\rho|)+S^{0}_{0,1}(|\rho|),$$
$$\partial_{x_1}^kA_N\in S^{2-3N-3k}_{2,2}(|\rho|)+S^{2}_{0,1}(|\rho|),\quad  
\partial_{x_1}^kB_N\in S^{3-4N-3k}_{2,2}(|\rho|)+S^{1-N}_{0,1}(|\rho|),$$ 
 with respect to the variables $x',\xi'$
uniformly in $z, h$ and $0\le x_1\le 2\delta\min\{1,\rho_1\}$. Finally, for $0<x_1\le 2\delta\min\{1,\rho_1\}$ we have
${\rm Im}\,\varphi\ge x_1{\rm Im}\,\rho/2$.
\end{prop}

Define the sets ${\cal M}_j\subset Z\times T^*\Gamma$, $j=1,2,$ 
as follows: ${\cal M}_1:=Z_{1,0}\times {\rm supp}\,\chi$, ${\cal M}_2:=Z_{1}\times {\rm supp}\,(1-\chi)\cup Z_2\times T^*\Gamma\cup
Z_3\times T^*\Gamma$. It follows from Lemma 3.1 that if $(z,x',\xi')\in{\cal M}_1$, then $C\sqrt{|{\rm Im}\,z|}\le|\rho|\le\widetilde C$ and 
${\rm Im}\,\rho\ge\frac{|{\rm Im}\,z|}{2|\rho|}$, while for $(z,x',\xi')\in{\cal M}_2$ we have 
$C_1\langle\xi'\rangle\le |\rho|\le C_2\langle\xi'\rangle$
and ${\rm Im}\,\rho\ge C\langle\xi'\rangle$.

Clearly, we have ${\cal D}_{x_1}\widetilde u_\psi|_{x_1=0}=T_\psi(z,h)f={\rm Op}_h(\tau_\psi)f$, where
$$\tau_\psi=a\frac{\partial\varphi}{\partial x_1}|_{x_1=0}-ih\frac{\partial a}{\partial x_1}|_{x_1=0}=\psi\rho-ih\sum_{j=0}^{N-1}
h^ja_{1,j}.$$

\begin{lemma} If $z\in Z_{1,0}$ we have the estimate
$$\left\|T_\psi(z,h)f-{\rm Op}_h(\psi\rho+hb_\psi)f\right\|_{H^1(\Gamma)}\le \frac{Ch}{\sqrt{|{\rm Im}\,z|}}
\|f\|_{L^2(\Gamma)}\eqno{(3.10)}$$
where $$b_\psi=-\frac{i}{2}(1-\chi)\psi q(0,x',1,\xi'/\sqrt{r_0(x',\xi')})-
\frac{1}{2}(1-\chi)\langle R(0,x')\xi'/\sqrt{r_0(x',\xi')},\nabla_{x'}\psi(x')\rangle.$$
Moreover, (3.10) holds for all $z\in Z_2\cup Z_3$ with $|{\rm Im}\,z|$ replaced by 1.
\end{lemma}

{\it Proof.} If $z\in Z_{1,0}$, it follows from the above
proposition that $\sum_{j=0}^{N-1}
h^j\chi a_{1,j}\in S^{-1}_{2,2}(\sqrt{|{\rm Im}\,z|})$, and hence by Proposition 2.1,
$$\left\|{\rm Op}_h(\sum_{j=0}^{N-1}h^j\chi a_{1,j})f\right\|_{H^1(\Gamma)}\le C
\left\|{\rm Op}_h(\sum_{j=0}^{N-1}h^j\chi a_{1,j})f\right\|_{L^2(\Gamma)}$$
$$\le \frac{C}{\sqrt{|{\rm Im}\,z|}}\left(1+\frac{\sqrt{h}}{|{\rm Im}\,z|}\right)^{d}\|f\|_{L^2(\Gamma)}\le  \frac{\widetilde C}{\sqrt{|{\rm Im}\,z|}}\|f\|_{L^2(\Gamma)}$$
as long as $|{\rm Im}\,z|\ge\sqrt{h}$. Clearly, the above bound holds for all $z\in Z_2\cup Z_3$ with $|{\rm Im}\,z|$ replaced by 1. 
On the other hand, it is easy to see that
$$(a_{1,0}-ib_\psi)(1-\chi)+\sum_{j=1}^{N-1}h^j(1-\chi)a_{1j}\in S^{-1}_{0,1}(\langle\xi'\rangle)={\cal S}_0^{-1}$$
 uniformly in $z$ and $h$. Hence the $h$-psdo with this symbol is bounded from $L^2$ to $H^1$.
\eproof

\begin{prop} Let $u_\psi$ satisfy $(P(h)-z)u_\psi=0$ in $\Omega$, $u_\psi|_{\Gamma}=\psi f$. Then, if $z\in Z_{1,0}$, 
$$\left\|\gamma{\cal D}_{\nu}u_\psi-T_\psi(z,h)f\right\|_{H^1(\Gamma)}\le C_Nh^{-s_d}\left(\frac{\sqrt{h}}{|{\rm Im}\,z|}\right)^{2N}\|f\|_{L^2(\Gamma)}\eqno{(3.11)}$$
with constants $C_N,s_d>0$ independent of $f$, $h$ and $z$, $s_d$ independent of $N$.
If $z\in Z_2\cup Z_3$, then (3.11) holds with $|{\rm Im}\,z|$ replaced by $1$.
\end{prop}

{\it Proof.} Given an integer $s\ge 0$, $H^s(\Omega)$ will denote the Sobolev space equipped with the semi-classical norm
$$\|g\|_{H^s(\Omega)}=\sum_{|\alpha|\le s}\|{\cal D}_x^\alpha g\|_{L^2(\Omega)}.$$
Denote also by $G_D$ the Dirichlet self-adjoint realization of the operator
$-n^{-1}\nabla c\nabla$ on the Hilbert space $L^2(\Omega,n(x)dx)$. Then the function
$$w_\psi:=u_\psi-\widetilde u_\psi+\left(h^2G_D-z\right)^{-1}\frac{c}{n}{\cal P}(z,h)\widetilde u_\psi$$
satisfies the equation $\left(h^2G_D-z\right)w_\psi=0$ in $\Omega$, $w_\psi|_{\Gamma}=0$. Since $z/h^2$ does not belong to the
spectrum of $G_D$, this implies that $w_\psi$ is identically zero. Thus we get
$$\left\|\gamma{\cal D}_\nu u_\psi-\gamma{\cal D}_\nu \widetilde u_\psi\right\|_{H^1(\Gamma)}\le 
\left\|\gamma{\cal D}_\nu\left(h^2G_D-z\right)^{-1}\frac{c}{n}{\cal P}(z,h)\widetilde u_\psi\right\|_{H^1(\Gamma)}$$
$$\le Ch^{-1/2}\left\|\left(h^2G_D-z\right)^{-1}\frac{c}{n}{\cal P}(z,h)\widetilde u_\psi\right\|_{H^4(\Omega)}\eqno{(3.12)}$$
where we have used the semi-classical version of the trace theorem. On the other hand, it is well known that the resolvent
of the operator $G_D$ satisfies the bound
$$\left\|\left(h^2G_D-z\right)^{-1}\right\|_{H^{2k}(\Omega)\to H^{2k}(\Omega)}\le \frac{C_k}{|{\rm Im}\,z|}\eqno{(3.13)}$$
for every integer $k\ge 0$. Indeed, for $k=0$ (3.13) is trivial, while for $k\ge 1$ it follows from the coercive estimate
$$\left\|v\right\|_{H^{2k}(\Omega)}\le \widetilde C_k\left\|h^2G_Dv\right\|_{H^{2k-2}(\Omega)}+\widetilde C_k\left\|v\right\|_{H^{2k-2}(\Omega)},\quad\forall
v\in D(G_D)\cap H^{2k-2}(\Omega).$$
Thus, (3.11) follows from (3.12), (3.13) and the following

\begin{prop} If $z\in Z_{1,0}$, given any integer $s\ge 0$ there are $\ell_s, N_s>0$ so that for $N\ge N_s$ we have the estimate
$$\left\|{\cal P}(z,h)\widetilde u_\psi\right\|_{H^s(\Omega)}\le C_Nh^{-\ell_s}\left(\frac{\sqrt{h}}{|{\rm Im}\,z|}\right)^{2N}\|f\|_{L^2(\Gamma)}.\eqno{(3.14)}$$ 
If $z\in Z_2\cup Z_3$, then (3.14) holds with $|{\rm Im}\,z|$ replaced by $1$.
\end{prop}

{\it Proof.} In view of (3.9) we can write
$${\cal P}(z,h)\widetilde u_\psi=(2\pi h)^{-d+1}\int\int e^{-\frac{i}{h}\langle x'-y',\xi'\rangle}
K(x,\xi',z,h)f(y')dy'd\xi',$$
where 
$$K=e^{\frac{i}{h}\langle x',\xi'\rangle}\left[{\cal P}(z,h),\Phi_\delta\right]e^{-\frac{i}{h}\langle x',\xi'\rangle}e^{\frac{i}{h}\widetilde\varphi}a
+e^{\frac{i}{h}\widetilde\varphi}\Phi_\delta\left(x_1^NA_N+h^NB_N\right)=:K_1+K_2.$$

\begin{lemma} If $z\in Z_{1,0}$, for any multi-index $\alpha$ there are $\ell_\alpha, N_\alpha>0$ so that for $N\ge N_\alpha$ we have
$$\left|\partial_{x}^{\alpha} K\right|\le C_{\alpha,N}h^{-\ell_\alpha}
\left(\frac{\sqrt{h}}{|{\rm Im}\,z|}\right)^{2N}.\eqno{(3.15)}$$
If $z\in Z_2\cup Z_3$, then (3.15) holds with $|{\rm Im}\,z|$ replaced by $1$.
\end{lemma}

{\it Proof.} An easy computation leads to the identity
$$\left[{\cal P}(z,h),\Phi_\delta\right]=-2ih\frac{\partial\Phi_\delta}{\partial x_1}{\cal D}_{x_1}
-2ih\left\langle R(x)\nabla_{x'}\Phi_\delta,{\cal D}_{x'}\right\rangle$$
$$-h^2\frac{\partial^2\Phi_\delta}{\partial x_1^2}-h^2\sum_{ij}R_{ij}(x)\frac{\partial^2\Phi_\delta}{\partial x'_i\partial x'_j}
-ih^2q(x,\nabla_{x}\Phi_\delta).$$
Hence
$$e^{\frac{i}{h}\langle x',\xi'\rangle}\left[{\cal P}(z,h),\Phi_\delta\right]e^{-\frac{i}{h}\langle x',\xi'\rangle}
=-2ih\frac{\partial\Phi_\delta}{\partial x_1}{\cal D}_{x_1}
-2ih\left\langle R(x)\nabla_{x'}\Phi_\delta,{\cal D}_{x'}\right\rangle$$
$$-h^2\frac{\partial^2\Phi_\delta}{\partial x_1^2}-h^2\sum_{ij} R_{ij}(x)\frac{\partial^2\Phi_\delta}{\partial x'_i\partial x'_j}
-ih^2q(x,\nabla_{x}\Phi_\delta)-2ih\left\langle R(x)\nabla_{x'}\Phi_\delta,\xi'\right\rangle.$$
Observe now that if $|\alpha|\ge 1$, then the function $\partial_{x}^{\alpha}\Phi_\delta$ is supported in the region
$\Theta:=\delta\min\{1,\rho_1\}\le x_1\le 2\delta\min\{1,\rho_1\}$.  To prove (3.15) we will consider two cases.

Case 1. $(z,x',\xi')\in{\cal M}_1$. Then $\rho_1=|\rho|^3$ and $C'|\rho|^3\le\min\{1,\rho_1\}\le|\rho|^3$.
It is easy to see that in this case we have $\partial_{x}^{\alpha}\Phi_\delta=O\left(|{\rm Im}\,z|^{-\ell_\alpha}\right)=
O\left(h^{-\ell_\alpha/2}\right)$ as long as $|{\rm Im}\,z|\ge \sqrt{h}$. In view of Proposition 3.4, on $\Theta$ we also have
$$\left|e^{\frac{i}{h}\widetilde\varphi}a\right|\le \widetilde C\exp\left(-{\rm Im}\,\widetilde\varphi/h\right)\le 
\widetilde C\exp\left(-\frac{x_1|{\rm Im}\,z|}{2h|\rho|}\right)$$ $$\le 
\widetilde C\exp\left(-C|\rho|^2|{\rm Im}\,z|/h\right)\le 
\widetilde C\exp\left(-C|{\rm Im}\,z|^2/h\right)$$
and more generally
$$\left|\partial_{x}^{\alpha}\left(e^{\frac{i}{h}\widetilde\varphi}a\right)\right|\le 
\widetilde C_\alpha h^{-\ell_\alpha}\exp\left(-C|{\rm Im}\,z|^2/h\right).$$
Thus we get
$$\left|\partial_{x}^{\alpha}K_1\right|\le 
\widetilde C_\alpha h^{-\ell_\alpha}\exp\left(-C|{\rm Im}\,z|^2/h\right)\eqno{(3.16)}$$
with probably new constants. Furthermore, on supp$\,\Phi_\delta$ we have
$$\left|x_1^Ne^{\frac{i}{h}\widetilde\varphi}\right|\le \widetilde Cx_1^N\exp\left(-\frac{x_1|{\rm Im}\,z|}{2h|\rho|}\right)\le 
C_N\left(\frac{h|\rho|}{|{\rm Im}\,z|}\right)^N.$$
On the other hand, by Proposition 3.4 we have $A_N=O_N\left(|\rho|^{-3N}\right)$, $B_N=O_N\left(|\rho|^{-4N}\right)$. Hence
$$\left|x_1^Ne^{\frac{i}{h}\widetilde\varphi}A_N\right|+\left|h^Ne^{\frac{i}{h}\widetilde\varphi}B_N\right| 
\le C_N\left(\frac{h}{|\rho|^2|{\rm Im}\,z|}\right)^N+C_N\left(\frac{h}{|\rho|^4}\right)^N\le
C_N\left(\frac{\sqrt{h}}{|{\rm Im}\,z|}\right)^{2N}.$$
Moreover, it is easy to see that differentiating these two functions makes appear additional factors 
$O\left(|{\rm Im}\,z|^{-\ell_1}h^{-\ell_2}\right)=O\left(h^{-\ell_1/2-\ell_2}\right)$ with $\ell_1$ and $\ell_2$
depending only on the order of differentiation. Thus we get
$$\left|\partial_{x}^{\alpha}K_2\right|\le C_N h^{-\ell_\alpha}\left(\frac{\sqrt{h}}{|{\rm Im}\,z|}\right)^{2N}.\eqno{(3.17)}$$
Clearly, in this case (3.15) follows from (3.16) and (3.17).

Case 2. $(z,x',\xi')\in{\cal M}_2$. Then $\min\{1,\rho_1\}\ge Const>0$. As above, it is easy to see that in this
case we have
$$\left|\partial_{x}^{\alpha}K_1\right|\le 
\widetilde C_\alpha\left(h^{-1}\langle\xi'\rangle\right)^{\ell_\alpha}\exp\left(-\frac{C\langle\xi'\rangle}{h}\right)\eqno{(3.18)}$$
and
$$\left|\partial_{x}^{\alpha}K_2\right|\le C_N\left(h^{-1}\langle\xi'\rangle\right)^{\ell_\alpha}\left(\frac{h}{\langle\xi'\rangle}\right)^N.\eqno{(3.19)}$$
If $N\ge \ell_\alpha$ we deduce from these bounds that $\partial_{x}^{\alpha}K_1, \partial_{x}^{\alpha}K_2=O_N\left(h^{N-\ell'_\alpha}\right)$, 
which again implies (3.15).

\eproof

It follows from Proposition 2.1 and Lemma 3.8 that, for $z\in Z_{1,0}$, 
$$\left\|\partial_{x}^{\alpha}{\cal P}(z,h)\widetilde u_\psi(x_1,\cdot)\right\|_{L^2(\Gamma)}\le C_Nh^{-\widetilde\ell_\alpha}\left(\frac{\sqrt{h}}{|{\rm Im}\,z|}\right)^{2N}\|f\|_{L^2(\Gamma)}$$ 
for $0\le x_1\le 2\delta$. Hence
$$\left\|\partial_{x}^{\alpha}{\cal P}(z,h)\widetilde u_\psi\right\|_{L^2(\Omega)}\le\left(\int_0^{2\delta}
\left\|\partial_{x}^{\alpha}{\cal P}(z,h)\widetilde u_\psi(x_1,\cdot)\right\|_{L^2(\Gamma)}^2dx_1\right)^{1/2}$$ $$\le \sqrt{2\delta} 
 C_Nh^{-\widetilde\ell_\alpha}\left(\frac{\sqrt{h}}{|{\rm Im}\,z|}\right)^{2N}\|f\|_{L^2(\Gamma)}$$ 
 which clearly implies (3.14) in this case. If $z\in Z_2\cup Z_3$, the above estimates clearly hold with $|{\rm Im}\,z|$ replaced by 1.
 \eproof

Let $\{\psi_j\}_{j=1}^J$ be a partition of the identity on $\Gamma$. Then we have
$u=\sum_{j=1}^Ju_{\psi_j}$ and $T(z,h)=\sum_{j=1}^J T_{\psi_j}(z,h)$ is an $h$-psdo on $\Gamma$ with
a principal symbol $\rho$. Observe that if $z\in Z_{1,\epsilon}$, there are $N_0=N_0(\epsilon)\gg 1$ and 
$h_0=h_0(\epsilon)\ll 1$ such that for $N\ge N_0$ and $0<h\le h_0$ we have
$$C_Nh^{-s_d}\left(\frac{\sqrt{h}}{|{\rm Im}\,z|}\right)^{2N}\le C_N h^{2\epsilon N-s_d}\le h.$$
This bound clearly holds for all $z\in Z_2\cup Z_3$ (with $|{\rm Im}\,z|$ replaced by 1). 
Therefore, (3.7) follows from (3.10) and (3.11) with $b=\sum_{j=1}^Jb_{\psi_j}$.
\eproof

In what follows, given any $s\in{\bf R}$ we denote $\|f\|_{H^s(\Gamma)}:=\|{\rm Op}_h(\langle\xi'\rangle^s)f\|_{L^2(\Gamma)}$.

\begin{lemma} Let $z\in Z_2$. Then we have 
$$\left\|\frac{dT}{dz}(z,h)f-{\rm Op}_h(\frac{d\rho}{dz}(z))f\right\|_{L^2(\Gamma)}\le Ch\|f\|_{H^{-1}(\Gamma)}
 \eqno{(3.20)}$$
 with a constant $C>0$ independent of $z$, $h$ and $f$. Moreover,
$$\left|{\rm Re}\,\left\langle cT(-1,h)f,f\right\rangle_{L^2(\Gamma)}\right|\le C_Nh^{N-s_d}\left\|f\right\|_{L^2(\Gamma)}^2.\eqno{(3.21)}$$
\end{lemma}

{\it Proof.} It follows from Lemma 4.3 below that
$$\sum_{j=0}^{N-1}h^j\frac{da_{1,j}}{dz}\in S^{-1}_{0,1}(\langle\xi'\rangle)={\cal S}_0^{-1}.$$
Hence the $h$-psdo with this symbol is bounded from $H^{-1}$ to $L^2$ uniformly in $h$ and $z\in Z_2$, which implies (3.20).
To prove (3.21) observe that by Green's formula we have the identity
$${\rm Im}\,\left\langle c\partial_\nu\widetilde u|_\Gamma,f\right\rangle_{L^2(\Gamma)}=-{\rm Im}\,\left\langle \nabla c\nabla \widetilde u,\widetilde u\right\rangle_{L^2(\Omega)}=-{\rm Im}\,\left\langle (\nabla c\nabla -h^{-2}n)\widetilde u,\widetilde u\right\rangle_{L^2(\Omega)}$$
where $\widetilde u=\sum_{j=1}^J\widetilde u_{\psi_j}$. Hence
$$\left|{\rm Re}\,\left\langle cT(-1,h)f,f\right\rangle_{L^2(\Gamma)}\right|\le h^{-1}\|{\cal P}(-1,h)\widetilde u\|_{L^2(\Omega)}
\|\widetilde u\|_{L^2(\Omega)}.\eqno{(3.22)}$$
By Proposition 2.1 it is easy to see that $\|\widetilde u\|_{L^2(\Omega)}\le Ch^{-s'_d}\|f\|_{L^2(\Gamma)}$, which together with (3.14)
and (3.22) imply (3.21).
\eproof

\section{Proof of Proposition 3.4} 
We will first solve equation (3.8). We can expand the functions $R(x)$ and $m(x)$ as follows
$$R(x)=\sum_{k=0}^{N-1}x_1^kR_k(x')+x_1^N{\cal R}_N(x),$$
$$m(x)=\sum_{k=0}^{N-1}x_1^km_k(x')+x_1^NM_N(x),$$
where $R_k$, ${\cal R}_N$, $m_k$, $M_N$ are smooth functions. Thus, if $\varphi=\sum_{k=0}^{N-1}x_1^k\varphi_k(x')$, we have
$${\cal E}:=(\partial_{x_1}\varphi)^2+\left\langle R(x)\nabla_{x'}\varphi,\nabla_{x'}\varphi\right\rangle -z m(x)$$ $$=
\sum_{k=0}^{N-2}\sum_{j=0}^{N-2}(k+1)(j+1)x_1^{k+j}\varphi_{k+1}\varphi_{j+1}$$ $$+\sum_{k=0}^{N-1}\sum_{j=0}^{N-1}x_1^{k+j}
\left\langle R\nabla_{x'}\varphi_k,\nabla_{x'}\varphi_j\right\rangle -z \sum_{k=0}^{N-1}x_1^km_k -zx_1^NM_N$$
$$=\sum_{k+j\le N-1}(k+1)(j+1)x_1^{k+j}\varphi_{k+1}\varphi_{j+1}$$ $$+\sum_{k+j\le N-1}x_1^{k+j}
\left\langle R\nabla_{x'}\varphi_k,\nabla_{x'}\varphi_j\right\rangle -z \sum_{k=0}^{N-1}x_1^km_k+x_1^N\Psi_N^{(1)}$$
where
$$\Psi_N^{(1)}=\sum_{k,j\le N-2,\, k+j\ge N}(k+1)(j+1)x_1^{k+j-N}\varphi_{k+1}\varphi_{j+1}$$
$$+\sum_{k,j\le N-1,\, k+j\ge N}x_1^{k+j-N}
\left\langle R\nabla_{x'}\varphi_k,\nabla_{x'}\varphi_j\right\rangle -z M_N.$$
We also have
$$\sum_{k+j\le N-1}x_1^{k+j}
\left\langle R\nabla_{x'}\varphi_k,\nabla_{x'}\varphi_j\right\rangle =\sum_{k+j+\ell\le N-1}x_1^{k+j+\ell}
\left\langle R_\ell\nabla_{x'}\varphi_k,\nabla_{x'}\varphi_j\right\rangle+ x_1^N\Psi_N^{(2)}$$
where
$$\Psi_N^{(2)}=\sum_{k+j\le N-1}x_1^{k+j}
\left\langle{\cal R}_N\nabla_{x'}\varphi_k,\nabla_{x'}\varphi_j\right\rangle$$ $$+
\sum_{\ell\le N-1,\,k+j\le N-1,\,k+j+\ell\ge N}x_1^{k+j+\ell-N}
\left\langle R_\ell\nabla_{x'}\varphi_k,\nabla_{x'}\varphi_j\right\rangle.$$
Thus we have ${\cal E}=x_1^N\Psi_N$ with $\Psi_N=\Psi_N^{(1)}+\Psi_N^{(2)}$, provided the coefficients $\varphi_k$ satisfy the 
relationships
$$\sum_{k+j=K}(k+1)(j+1)\varphi_{k+1}\varphi_{j+1}+\sum_{k+j+\ell=K}
\left\langle R_\ell\nabla_{x'}\varphi_k,\nabla_{x'}\varphi_j\right\rangle -zm_K=0\eqno{(4.1)}$$
for every integer $0\le K\le N-2$. Clearly, if we take $\varphi_0=-\langle x'-y',\xi'\rangle$, then $\varphi_1=\rho$ is a
solution of (4.1) with $K=0$. Now, given $\varphi_j$, $0\le j\le K-1$, $K\ge 2$, we can determine $\varphi_K$ in a unique way by (4.1).

\begin{lemma} We have $\varphi_k\in S^{4-3k}_{2,2}(|\rho|)+S^{1}_{0,1}(|\rho|)$, $1\le k\le N-1$, 
$\partial_{x_1}^k\Psi_N\in S^{2-3N-3k}_{2,2}(|\rho|)+S^{2}_{0,1}(|\rho|)$, $k\ge 0$, uniformly in $z$ and $0\le x_1\le 2\delta\min\{1,|\rho|^3\}$. Moreover, if $\delta>0$ is small enough, independent of $\rho$, we have 
$${\rm Im}\,\varphi\ge x_1{\rm Im}\,\rho/2\quad \mbox{for}\quad 0\le x_1\le 2\delta\min\{1,|\rho|^3\}.\eqno{(4.2)}$$
\end{lemma}

{\it Proof.} In view of Lemma 3.2
we have $zm_K\rho^{-1}\in S^{-1}_{2,2}(|\rho|)+S^{-1}_{0,1}(|\rho|)$ uniformly in $z$. We will now proceed by induction. 
Suppose that $\varphi_k\in S^{4-3k}_{2,2}(|\rho|)+S^{1}_{0,1}(|\rho|)$, $1\le k\le K$.
This implies $\nabla_{x'}\varphi_k\in S^{2-3k}_{2,2}(|\rho|)+S^{1}_{0,1}(|\rho|)$, $1\le k\le K$,
which yields 
$$\left\langle R_\ell\nabla_{x'}\varphi_k,\nabla_{x'}\varphi_j\right\rangle\in S^{4-3K}_{2,2}(|\rho|)+S^{2}_{0,1}(|\rho|),\quad 
k+j+\ell\le K,\,k,j\ge 1.\eqno{(4.3)}$$
Furthermore, since $\nabla_{x'}\varphi_0=\xi'$, we have 
$$\left\langle R_\ell\nabla_{x'}\varphi_k,\nabla_{x'}\varphi_0\right\rangle,\,
\left\langle R_\ell\nabla_{x'}\varphi_0,\nabla_{x'}\varphi_j\right\rangle\in
S^{2-3K}_{2,2}(|\rho|)+S^{2}_{0,1}(|\rho|),\quad 
1\le k,j\le K,\eqno{(4.4)}$$
$$\left\langle R_\ell\nabla_{x'}\varphi_0,\nabla_{x'}\varphi_0\right\rangle\in
S^{0}_{2,2}(|\rho|)+S^{2}_{0,1}(|\rho|).\eqno{(4.5)}$$
We also have 
$$\varphi_{k+1}\varphi_{j+1}\in S^{2-3K}_{2,2}(|\rho|)+S^{2}_{0,1}(|\rho|),\quad 
k+j= K,\,k,j\ge 1.\eqno{(4.6)}$$
Thus by equation (4.1) and (4.3)-(4.6) we conclude that $2\rho\varphi_{K+1}-zm_K\in 
S^{2-3K}_{2,2}(|\rho|)+S^{2}_{0,1}(|\rho|)$, and hence
$\varphi_{K+1}\in S^{1-3K}_{2,2}(|\rho|)+S^{1}_{0,1}(|\rho|)$ as desired.
The property concerning the function $\Psi_N$ follows easily from the following observation: if $k\ge 0$ and
$a\in S^{\ell_1}_{2,2}(|\rho|)+S^{\ell_2}_{0,1}(|\rho|)$, then for $0\le x_1\le 2\delta\min\{1,|\rho|^3\}$ we have
$x_1^ka\in S^{\ell_1+3k}_{2,2}(|\rho|)+S^{\ell_2}_{0,1}(|\rho|)$.

To bound ${\rm Im}\,\varphi$ from below we will show that for every multi-index $\alpha$ we have the estimate
$$\left|{\rm Im}\,\partial_{x'}^\alpha\varphi_k\right|\le \frac{C_{k,\alpha}{\rm Im}\,\rho}{\min\{1,|\rho|^{3k-3+2|\alpha|}\}},
\quad k\ge 1.\eqno{(4.7)}$$
For $(z,x',\xi')\in{\cal M}_2$ we have 
$$\left|{\rm Im}\,\partial_{x'}^\alpha\varphi_k\right|\le\left|\partial_{x'}^\alpha\varphi_k\right|\le 
\widetilde C_{k,\alpha}|\rho|\le C_{k,\alpha}{\rm Im}\,\rho$$
which implies (4.7) in this case. Let now $(z,x',\xi')\in{\cal M}_1$. Observe first that
$$|{\rm Im}\,(z\rho^{-1})|\le |{\rm Im}\,z||\rho|^{-1}+|z||{\rm Im}\,(\rho^{-1})|\le 2{\rm Im}\,\rho+2|\rho|^{-2}{\rm Im}\,\rho
\le C|\rho|^{-2}{\rm Im}\,\rho.$$
Differentiating equation (4.1) we obtain
$$2(K+1)\partial_{x'}^\alpha\varphi_{K+1}+\rho^{-1}\sum_{\Theta_{K,\alpha}}(k+1)(j+1)\partial_{x'}^{\alpha_1}\varphi_{k+1}\partial_{x'}^{\alpha_2}\varphi_{j+1}$$ $$+\rho^{-1}\sum_{\widetilde\Theta_{K,\alpha}}
\left\langle \partial_{x'}^{\alpha_1}R_\ell\partial_{x'}^{\alpha_2}
\nabla_{x'}\varphi_k,\partial_{x'}^{\alpha_3}\nabla_{x'}\varphi_j\right\rangle -z\rho^{-1}\partial_{x'}^{\alpha}m_K=0\eqno{(4.8)}$$
where $\Theta_{K,\alpha}:=\{(k,j,\alpha_1,\alpha_2):k+j=K,k,j\le K-1, |\alpha_1|+|\alpha_2|=|\alpha|;|\alpha_1|, |\alpha_2|\le |\alpha|-1\}$,
$\widetilde\Theta_{K,\alpha}:=\{(\ell,k,j,\alpha_1,\alpha_2,\alpha_3):\ell+k+j=K, |\alpha_1|+|\alpha_2|+|\alpha_3|=|\alpha|\}$.
To prove (4.7) in this case we will proceed by induction in $k$ and $|\alpha|$. Fix integers $K\ge 1$ and $A\ge 0$ and suppose that
(4.7) holds for $1\le k\le K$ and all $\alpha$, and for $k=K+1$ and $|\alpha|\le A-1$. We have to show that (4.7) holds for $k=K$
and $|\alpha|=A$. To this end we will use (4.8). Observe that on $\Theta_{K,\alpha}$ we have
$$\left|{\rm Im}\,\left(\rho^{-1}\partial_{x'}^{\alpha_1}\varphi_{k+1}\partial_{x'}^{\alpha_2}\varphi_{j+1}\right)\right|\le
\left|{\rm Im}\,\left(\rho^{-1}\right)\right|\left|\partial_{x'}^{\alpha_1}\varphi_{k+1}\right|
\left|\partial_{x'}^{\alpha_2}\varphi_{j+1}\right|$$ $$+|\rho|^{-1}\left|{\rm Im}\,\partial_{x'}^{\alpha_1}\varphi_{k+1}\right|\left|\partial_{x'}^{\alpha_2}\varphi_{j+1}\right|+|\rho|^{-1}\left|\partial_{x'}^{\alpha_1}\varphi_{k+1}\right|\left|{\rm Im}\,\partial_{x'}^{\alpha_2}\varphi_{j+1}\right|\le C_{K,\alpha}|\rho|^{-3K-2A}{\rm Im}\,\rho\eqno{(4.9)}$$
where we have used our hypothesis and the fact that in this case $\varphi_k\in S_{2,2}^{4-3k}(|\rho|)$.
Similarly, on $\widetilde\Theta_{K,\alpha}$, $k,j\ge 1$, we have
$$\left|{\rm Im}\,\left(\rho^{-1}\partial_{x'}^{\alpha_2}
\partial_{x'}\varphi_k\partial_{x'}^{\alpha_3}\partial_{x'}\varphi_j\right)\right|\le C_{K,\alpha}|\rho|^{2-3K-2A}{\rm Im}\,\rho.\eqno{(4.10)}$$
If one of $k$ or $j$ is 0, since $\nabla_{x'}\varphi_0=\xi'$ is bounded on supp$\,\chi$, the left-hand side of (4.10) is  
$O\left(|\rho|^{-3K-2A}{\rm Im}\,\rho\right)$, while for $j=k=0$ it is $O\left({\rm Im}\,(\rho^{-1})\right)=O\left(|\rho|^{-2}{\rm Im}\,\rho\right)$.
Thus, by (4.8) we conclude that
$$\left|{\rm Im}\,\partial_{x'}^{\alpha}\varphi_{K+1}\right|\le C_{K,\alpha}|\rho|^{-3K-2A}{\rm Im}\,\rho$$
which is the desired bound.

Using (4.7) with $\alpha=0$ we obtain, for $0<x_1\le 2\delta\min\{1,|\rho|^3\}$, 
$${\rm Im}\,\varphi\ge x_1{\rm Im}\,\rho-x_1\sum_{k=2}^{N-1}x_1^{k-1}|{\rm Im}\,\varphi_k|
\ge x_1{\rm Im}\,\rho-x_1{\rm Im}\,\rho\sum_{k=2}^{N-1}C_kx_1^{k-1}(\min\{1,|\rho|^3\})^{-k+1}$$
 $$\ge x_1{\rm Im}\,\rho\left(1-O(\delta)\right)\ge x_1{\rm Im}\,\rho/2$$
 provided $\delta>0$ is taken small enough, independent of $\rho$.
\eproof

To solve equation (3.9) observe first that
$${\cal Q}:=e^{-\frac{i}{h}\varphi}{\cal P}(z,h)e^{\frac{i}{h}\varphi}={\cal P}(z,h)$$ $$+2\frac{\partial\varphi}{\partial x_1}{\cal D}_{x_1}
+2\left\langle R(x)\nabla_{x'}\varphi,{\cal D}_{x'}\right\rangle+hq(x,\nabla_x\varphi)
+\left(\frac{\partial\varphi}{\partial x_1}\right)^2+r(x,\nabla_{x'}\varphi)$$
 $$={\cal D}_{x_1}^2+r(x,{\cal D}_{x'})+hq(x,{\cal D}_x)+h^2\widetilde q(x)$$
 $$+2\frac{\partial\varphi}{\partial x_1}{\cal D}_{x_1}
+2\left\langle R(x)\nabla_{x'}\varphi,{\cal D}_{x'}\right\rangle+hq(x,\nabla_x\varphi)
+x_1^N\Psi_N$$
where we have used that the phase function satisfies equation (3.8). Write
$$q(x,\xi)=\sum_{k=0}^{N-1}x_1^kq_k(x',\xi)+x_1^NQ_N(x,\xi),$$
$$q_k(x',\xi)=q_k^\sharp(x')\xi_1+q_k^\flat(x',\xi'),$$
$$\widetilde q(x)=\sum_{k=0}^{N-1}x_1^k\widetilde q_k(x')+x_1^N\widetilde Q_N(x).$$
We will be searching a solution to (3.9) in the form
$a=\sum_{j=0}^{N-1}h^ja_j(x,z)$, $a_0|_{x_1=0}=\psi$, $a_j|_{x_1=0}=0$, $j\ge 1$. Thus, if the functions $a_j$ satisfy the transport equations
$$-2i\frac{\partial\varphi}{\partial x_1}\frac{\partial a_j}{\partial x_1}
-2i\left\langle R(x)\nabla_{x'}\varphi,\nabla_{x'}a_j\right\rangle+q(x,\nabla_x\varphi)a_j$$
 $$=\left(\partial_{x_1}^2+r(x,\partial_{x'})+iq(x,\partial_x)-\widetilde q(x)\right)a_{j-1}+x_1^NA_N^{(j)},\quad 0\le j\le N-1,\eqno{(4.11)}$$
 $a_{-1}=0$, then 
$${\cal Q}a=x_1^N\sum_{j=0}^{N}h^jA_N^{(j-1)}-h^N\left(\partial_{x_1}^2+r(x,\partial_{x'})+iq(x,\partial_x)-\widetilde q(x)\right)a_{N-1}=
x_1^NA_N+h^NB_N
\eqno{(4.12)}$$
where we have put $A_N^{(-1)}=\Psi_N$. We will be looking for solutions of (4.11) in the form $a_j=\sum_{k=0}^{N-1}x_1^ka_{k,j}$,
$a_{0,0}=\psi$, $a_{0,j}=0$, $j\ge 1$. We have
$$\frac{\partial\varphi}{\partial x_1}\frac{\partial a_j}{\partial x_1}=\sum_{\nu+k\le N-1} x_1^{\nu+k}(\nu+1)(k+1)
\varphi_{\nu+1}a_{k+1,j}$$ $$+x_1^N\sum_{\nu+k\ge N, \nu,k\le N-1} x_1^{\nu+k-N}(\nu+1)(k+1)
\varphi_{\nu+1}a_{k+1,j},$$
$$\left\langle R(x)\nabla_{x'}\varphi,\nabla_{x'}a_j\right\rangle=\sum_{\ell+\nu+k\le N-1}x_1^{\ell+\nu+k}
\left\langle R_\ell(x')\nabla_{x'}\varphi_\nu,\nabla_{x'}a_{k,j}\right\rangle$$
 $$+x_1^N\sum_{\ell+\nu+k\ge N}x_1^{\ell+\nu+k-N}
\left\langle R_\ell(x')\nabla_{x'}\varphi_\nu,\nabla_{x'}a_{k,j}\right\rangle+x_1^N\left\langle{\cal R}_N(x)\nabla_{x'}\varphi,\nabla_{x'}a_j\right\rangle,$$
 $$q(x,\nabla_x\varphi)a_j=\sum_{\ell=0}^{N-1}x_1^\ell q_\ell(x',\nabla_x\varphi)a_j+x_1^NQ_N(x,\nabla_x\varphi)a_j$$
  $$=\sum_{\ell=0}^{N-1}x_1^\ell\left(q_\ell^\sharp(x')\partial_{x_1}\varphi
  +q^\flat_\ell(x',\nabla_{x'}\varphi)\right)a_j+x_1^NQ_N(x,\nabla_x\varphi)a_j$$
   $$=\sum_{\ell+\nu+k\le N-1}x_1^{\ell+\nu+k}\left((\nu+1)q_\ell^\sharp(x')\varphi_{\nu+1}+q^\flat_\ell(x',\nabla_{x'}\varphi_\nu)\right)a_{k,j}$$
  $$+x_1^N\sum_{\ell+\nu+k\ge N}x_1^{\ell+\nu+k-N}\left((\nu+1)q_\ell^\sharp(x')\varphi_{\nu+1}+
  q^\flat_\ell(x',\nabla_{x'}\varphi_\nu)\right)a_{k,j}+x_1^NQ_N(x,\nabla_x\varphi)a_j,$$
  $$\left(\partial_{x_1}^2+r(x,\partial_{x'})+iq(x,\partial_x)-\widetilde q(x)\right)a_{j-1}$$ $$=\sum_{k=0}^{N-2}x_1^k(k+2)(k+1)a_{k+2,j-1}
  +\sum_{\ell+k\le N-1}x_1^{\ell+k}\left\langle R_\ell(x')\nabla_{x'},\nabla_{x'}a_{k,j-1}\right\rangle$$ $$
  +\sum_{\ell+k\le N-1}x_1^{k+\ell}\left((k+1)iq_\ell^\sharp a_{k+1,j-1}+iq^\flat_\ell(x',\nabla_{x'}a_{k,j-1})-\widetilde q_\ell a_{k,j-1}\right)$$
  $$+x_1^N\sum_{\ell+k\ge N}x_1^{\ell+k-N}\left\langle R_\ell(x')\nabla_{x'},\nabla_{x'}a_{k,j-1}\right\rangle$$ $$
  +x_1^N\sum_{\ell+k\ge N}x_1^{k+\ell-N}\left((k+1)iq_\ell^\sharp a_{k+1,j-1}+iq^\flat_\ell(x',\nabla_{x'}a_{k,j-1})-\widetilde q_\ell a_{k,j-1}\right)$$ $$+x_1^N\left\langle{\cal R}_N(x')\nabla_{x'},\nabla_{x'}a_{j-1}\right\rangle+x_1^N\left(
  iQ_N(x,\nabla_{x}a_{j-1})-\widetilde Q_Na_{j-1}\right).$$
  Thus we obtain that the coefficients $a_{k,j}$ must satisfy the equations
  $$-2i\sum_{\nu+k=K}(\nu+1)(k+1)\varphi_{\nu+1}a_{k+1,j}-2i\sum_{\ell+\nu+k=K}
\left\langle R_\ell(x')\nabla_{x'}\varphi_\nu,\nabla_{x'}a_{k,j}\right\rangle$$
  $$+\sum_{\ell+\nu+k=K}\left((\nu+1)q_\ell^\sharp(x')\varphi_{\nu+1}+q^\flat_\ell(x',\nabla_{x'}\varphi_\nu)\right)a_{k,j}$$
  $$=(K+2)(K+1)a_{K+2,j-1}
  +\sum_{\ell+k=K}\left\langle R_\ell(x')\nabla_{x'},\nabla_{x'}a_{k,j-1}\right\rangle$$ $$
  +\sum_{\ell+k=K}\left((k+1)iq_\ell^\sharp a_{k+1,j-1}+iq^\flat_\ell(x',\nabla_{x'}a_{k,j-1})-\widetilde q_\ell a_{k,j-1}\right)\eqno{(4.13)}$$
  for every integer $0\le K\le N-2$. Clearly, there exist unique solutions $a_{k,j}$ of (4.13) such that $a_{0,0}=\psi$, 
  $a_{0,j}=0$, $j\ge 1$, and $a_{k,-1}=0$, $k\ge 0$. 
  
  \begin{lemma} We have $a_{k,j}\in S^{2-3k-4j}_{2,2}(|\rho|)+S^{-j}_{0,1}(|\rho|)$, $k\ge 1$, $j\ge 0$, 
  $\partial_{x_1}^kA_N^{(j)}\in S^{-3N-4j-3k}_{2,2}(|\rho|)+S^{1}_{0,1}(|\rho|)$, $j\ge 0$, 
  $\partial_{x_1}^kA_N\in S^{2-3N-3k}_{2,2}(|\rho|)+S^{2}_{0,1}(|\rho|)$, $\partial_{x_1}^kB_N\in S^{3-4N-3k}_{2,2}(|\rho|)+S^{1-N}_{0,1}(|\rho|)$,
  $k\ge 0$.
  \end{lemma}
  
  {\it Proof.} Observe first that equation (4.13) with $K=0$, $j=0$ yields the formula
  $$a_{1,0}=-\frac{i}{2} q(0,x',1,\xi'/\rho)\psi-\frac{1}{2\rho}\langle R(0,x')\xi',\nabla_{x'}\psi(x')\rangle.$$
  Hence $a_{1,0}\in S^{-1}_{2,2}(|\rho|)+S^{0}_{0,1}(|\rho|)$. To prove the assertion concerning the functions $a_{k,j}$ we will proceed by induction in $k$ and $j$. Fix $K\ge 1$, $J\ge 1$ and suppose that our assertion is true for all $0\le j\le J-1$, $k\ge 1$, and for $j=J$ and $1\le k\le K$.
  We have to show that it is true for $j=J$ and $k=K+1$. Our hypothesis together with Lemma 4.1 imply
  $$\varphi_{\nu+1}a_{k+1,J}\in S^{-3K-4J}_{2,2}(|\rho|)+S^{1-J}_{0,1}(|\rho|),\quad \nu+k=K,\,0\le k\le K-1,$$
  $$\left\langle R_\ell(x')\nabla_{x'}\varphi_\nu,\nabla_{x'}a_{k,J}\right\rangle\in S^{2-3K-4J}_{2,2}(|\rho|)+S^{1-J}_{0,1}(|\rho|),\quad
  \ell+\nu+k=K,\,\nu\ge 1,$$
  $$\left\langle R_\ell(x')\nabla_{x'}\varphi_0,\nabla_{x'}a_{k,J}\right\rangle\in S^{-3K-4J}_{2,2}(|\rho|)+S^{1-J}_{0,1}(|\rho|),\quad
  \ell+k=K,$$
  $$q_\ell^\sharp(x')\varphi_{\nu+1}a_{k,J}\in S^{3-3K-4J}_{2,2}(|\rho|)+S^{1-J}_{0,1}(|\rho|),\quad \ell+\nu+k=K,$$
  $$q^\flat_\ell(x',\nabla_{x'}\varphi_\nu)a_{k,J}\in S^{4-3K-4J}_{2,2}(|\rho|)+S^{1-J}_{0,1}(|\rho|),\quad
  \ell+\nu+k=K,\,\nu\ge 1,$$
  $$q^\flat_\ell(x',\nabla_{x'}\varphi_0)a_{k,J}\in S^{2-3K-4J}_{2,2}(|\rho|)+S^{1-J}_{0,1}(|\rho|),\quad
  \ell+k=K.$$
  One can also easily see that the right-hand side of equation (4.13) belongs to $S^{-3K-4J}_{2,2}(|\rho|)+S^{1-J}_{0,1}(|\rho|)$.
  Thus, by (4.13) we conclude that $\rho a_{K+1,J}\in S^{-3K-4J}_{2,2}(|\rho|)+S^{1-J}_{0,1}(|\rho|)$, which implies
  $a_{K+1,J}\in S^{-1-3K-4J}_{2,2}(|\rho|)+S^{-J}_{0,1}(|\rho|)$, as desired. 
  The properties concerning the functions $A_N^{(j)}$, $A_N$, $B_N$ follow easily from the following observation: if $k\ge 0$, $j\ge 0$, and
$a\in S^{\ell_1}_{2,2}(|\rho|)+S^{\ell_2}_{0,1}(|\rho|)$, then for $0\le x_1\le 2\delta\min\{1,|\rho|^3\}$ we have
$h^jx_1^ka\in S^{\ell_1+3k+4j}_{2,2}(|\rho|)+S^{\ell_2}_{0,1}(|\rho|)$, where we have used that $h\le |{\rm Im}\,z|^2\le C|\rho|^4$.
  \eproof
  
  \begin{lemma} Let $z\in Z_2$. Then $\frac{d\varphi_k}{dz}, \frac{da_{k,j}}{dz}\in S^{-1}_{0,1}(\langle\xi'\rangle)$, $k\ge 1$, $j\ge 0$.
  \end{lemma}
  
  {\it Proof.} Recall that in this case we have $C_1\langle\xi'\rangle\le |\rho|\le C_2\langle\xi'\rangle$. We also have $\frac{d\varphi_0}{dz}=0$
  and $2\rho\frac{d\rho}{dz}=-m_0(x')$. Hence $\frac{d\rho}{dz}\in S^{-1}_{0,1}(\langle\xi'\rangle)$. 
  Differentiating
  equation (4.1) once with respect to the variable $z$ it is easy to see that $\rho\frac{d\varphi_{K+1}}{dz}\in S^{0}_{0,1}(\langle\xi'\rangle)$,
  provided $\frac{d\varphi_k}{dz}\in S^{-1}_{0,1}(\langle\xi'\rangle)$ for $1\le k\le K$, which implies 
  $\frac{d\varphi_{K+1}}{dz}\in S^{-1}_{0,1}(\langle\xi'\rangle)$. Thus we obtain the desired properties of the functions $\frac{d\varphi_k}{dz}$
  by induction in $k$. Similalry, we have $\frac{da_{0,j}}{dz}=0$, $j\ge 0$, and $\frac{da_{1,0}}{dz}\in S^{-1}_{0,1}(\langle\xi'\rangle)$.
  Differentiating
  equation (4.13) once with respect to the variable $z$ it is easy to see that $\rho\frac{da_{K+1,J}}{dz}\in S^{0}_{0,1}(\langle\xi'\rangle)$,
  provided $\frac{da_{k,j}}{dz}\in S^{-1}_{0,1}(\langle\xi'\rangle)$ for $0\le j\le J-1$, $k\ge 1$, and $j=J$, $1\le k\le K$. 
  Therefore, the desired result follows by induction in $j$ and $k$.
  \eproof
  
\section{Eigenvalue-free regions}

In this section we will study the problem
$$\left\{
\begin{array}{lll}
\left(P_1(h)-z\right)u_1=0 &\mbox{in} &\Omega,\\
\left(P_2(h)-z\right)u_2=0 &\mbox{in} &\Omega,\\
u_1=u_2,\,\,\, c_1\partial_\nu u_1=c_2\partial_\nu u_2& \mbox{on}& \Gamma,
\end{array}
\right.
\eqno{(5.1)}
$$
where $z\in Z$, $0<h\ll 1$, $P_j(h)$, $j=1,2$, is defined by replacing in the definition of the operator
$P(h)$ from Section 3 the pair $(c,n)$ by $(c_j,n_j)$. Similarly, we define the functions $\rho_j$ by 
replacing in the definition of $\rho$ the function $m$ by $m_j=\frac{n_j}{c_j}$. We will also use the function $\chi$ introduced at the begining
of Section 3. Note that we can make the support of $\chi$ as large as we want by taking the parameter $\delta_0$ small enough. 
It follows from Theorem 3.3 that if
$z\in Z_{1,\epsilon}$, 
the function $f:=u_1|_\Gamma=u_2|_\Gamma$ satisfies the estimate
$$\left\|{\rm Op}_h(c_1\rho_1-c_2\rho_2)f\right\|_{L^2(\Gamma)}\le \frac{Ch}{\sqrt{|{\rm Im}\,z|}}
\|f\|_{L^2(\Gamma)}\eqno{(5.2)}$$
while in the case $c_1|_\Gamma\equiv c_2|_\Gamma$ we have the better estimate
$$\left\|{\rm Op}_h(\rho_1-\rho_2)f\right\|_{H^1(\Gamma)}\le \frac{Ch}{\sqrt{|{\rm Im}\,z|}}
\|f\|_{L^2(\Gamma)}.\eqno{(5.3)}$$
Moreover, (5.2) and (5.3) hold for all $z\in Z_2\cup Z_3$ with $|{\rm Im}\,z|$ replaced by 1.
We would like to invert the operators in the left-hand sides of (5.2) and (5.3). Note that it follows from Lemmas 3.1 and 3.2 that the function
$\rho_j$ satisfies the bounds, for $(z,x',\xi')\in{\cal M}_1$,
$$\left|\partial_{x'}^\alpha\partial_{\xi'}^\beta\rho_j\right|\le C_{\alpha,\beta}|{\rm Im}\,z|^{\frac{1}{2}-|\alpha|-|\beta|},\quad
|\alpha|+|\beta|\ge 1,\eqno{(5.4)}$$
$|\rho_j|\le Const$, while for $(z,x',\xi')\in{\cal M}_2$ we have 
$$\left|\partial_{x'}^\alpha\partial_{\xi'}^\beta\rho_j\right|\le C_{\alpha,\beta}\langle\xi'\rangle^{1-|\beta|}.\eqno{(5.5)}$$
 In particular, these estimates imply that $\rho_j\in {\cal S}^1_{\frac{1}{2}-\epsilon}$ if 
$z\in Z_{1,\epsilon}$ and $\rho_j\in {\cal S}^1_0$ if $z\in Z_2\cup Z_3$. Observe now that
$$c_1\rho_1-c_2\rho_2=\frac{\widetilde c(x')(c_0(x')r_0(x',\xi')-z)}{c_1\rho_1+c_2\rho_2}\eqno{(5.6)}$$
where $\widetilde c$ and $c_0$ are the restrictions on $\Gamma$ of the functions
$$c_1n_1-c_2n_2\quad\mbox{and}\quad\frac{c_1^2- c_2^2}{c_1n_1-c_2n_2}$$
respectively. It follows easily from (5.4)-(5.6) that
$$\left|\partial_{x'}^\alpha\partial_{\xi'}^\beta(c_1\rho_1-c_2\rho_2)\right|\le C_{\alpha,\beta}|{\rm Im}\,z|^{\frac{1}{2}-|\alpha|-|\beta|},\quad |\alpha|+|\beta|\ge 1,\eqno{(5.7)}$$
for $(z,x',\xi')\in{\cal M}_1$, and 
$$\left|\partial_{x'}^\alpha\partial_{\xi'}^\beta(c_1\rho_1-c_2\rho_2)\right|\le C_{\alpha,\beta}\langle\xi'\rangle^{k-|\beta|}\eqno{(5.8)}$$
for $(z,x',\xi')\in{\cal M}_2$ and all multi-indices $\alpha$ and $\beta$, where $k=-1$ if $c_0\equiv 0$,
$k=1$ if $c_0(x')\neq 0$, $\forall x'\in \Gamma$. In particular, these estimates imply that $c_1\rho_1-c_2\rho_2\in {\cal S}^k_{\frac{1}{2}-\epsilon}$ if 
$z\in Z_{1,\epsilon}$ and $c_1\rho_1-c_2\rho_2\in {\cal S}^k_0$ if $z\in Z_2\cup Z_3$. We will now consider two cases.

Case 1. $c_0\equiv 0$. Then $k=-1$. In this case we have $|\rho_1-\rho_2|\ge C\langle\xi'\rangle^{-1}$, $C>0$, so
$(\rho_1-\rho_2)^{-1}\in {\cal S}^1_{\frac{1}{2}-\epsilon}$ if 
$z\in Z_{1,\epsilon}$ and $(\rho_1-\rho_2)^{-1}\in {\cal S}^1_0$ if $z\in Z_2\cup Z_3$. Hence
$$\left\|{\rm Op}_h\left((\rho_1-\rho_2)^{-1}\right)g\right\|_{L^2(\Gamma)}\le C
\|g\|_{H^1(\Gamma)},\quad\forall g\in H^1(\Gamma),\eqno{(5.9)}$$
for $z\in Z_{1,\epsilon}\cup Z_2\cup Z_3$. By (5.3) and (5.9), for $z\in Z_{1,\epsilon}$, 
$$\left\|{\rm Op}_h\left((\rho_1-\rho_2)^{-1}\right){\rm Op}_h(\rho_1-\rho_2)f\right\|_{L^2(\Gamma)}\le \frac{Ch}{\sqrt{|{\rm Im}\,z|}}
\|f\|_{L^2(\Gamma)}.\eqno{(5.10)}$$
For $z\in Z_2\cup Z_3$, (5.10) holds with $|{\rm Im}\,z|$ replaced by 1. 
On the other hand, by Proposition 2.2 we have
$$\left\|{\rm Op}_h\left((\rho_1-\rho_2)^{-1}\right){\rm Op}_h(\rho_1-\rho_2)-Id\right\|_{L^2(\Gamma)\to L^2(\Gamma)}\le Ch^{2\epsilon}
.\eqno{(5.11)}$$
Combining (5.10) and (5.11) we conclude $\|f\|_{L^2}\le O(h^{2\epsilon})\|f\|_{L^2}$ for $z\in Z_{1,\epsilon}\cup Z_2\cup Z_3$,
which implies $f\equiv 0$ provided $h$ is taken small enough.

Case 2. $c_0(x')\neq 0$, $\forall x'\in \Gamma$. Then $k=1$. Observe first that the condition (1.7) implies $c_0>0$. It is easy to see that if  $z\in Z_2$, then 
we have $|c_1\rho_1-c_2\rho_2|\ge C\langle\xi'\rangle$, $C>0$, so  
 $(c_1\rho_1-c_2\rho_2)^{-1}\in {\cal S}^{-1}_0$. Hence  
$$\left\|{\rm Op}_h\left((c_1\rho_1-c_2\rho_2)^{-1}\right)g\right\|_{L^2(\Gamma)}\le C
\|g\|_{L^2(\Gamma)},\quad\forall g\in L^2(\Gamma).\eqno{(5.12)}$$
 By (5.2) and (5.12), if $z\in Z_2$, 
$$\left\|{\rm Op}_h\left((c_1\rho_1-c_2\rho_2)^{-1}\right){\rm Op}_h(c_1\rho_1-c_2\rho_2)f\right\|_{L^2(\Gamma)}\le Ch
\|f\|_{L^2(\Gamma)}.\eqno{(5.13)}$$
On the other hand, by Proposition 2.2 we have
$$\left\|{\rm Op}_h\left((c_1\rho_1-c_2\rho_2)^{-1}\right){\rm Op}_h(c_1\rho_1-c_2\rho_2)-Id\right\|_{L^2(\Gamma)\to L^2(\Gamma)}\le Ch^{2\epsilon}
.\eqno{(5.14)}$$
In the same way as above one can derive from (5.13) and (5.14) that $f\equiv 0$, provided $h$ is taken small enough.

Under the conditions (1.2) and (1.4) only, we have $|c_0r_0|\ge C|\xi'|^2$, $C>0$. Given any $0<\delta'\ll 1$ and any multi-indices $\alpha$ and $\beta$, by induction in $|\alpha|+|\beta|$ one can easily prove that the following estimates hold true:
$$\left|\partial_{x'}^\alpha\partial_{\xi'}^\beta\left((c_0r_0-z)^{-1}\right)\right|\le C_{\alpha,\beta}|{\rm Im}\,z|^{-1-|\alpha|-|\beta|}\eqno{(5.15)}$$
for $|c_0r_0-{\rm Re}\,z|\le\delta'$, ${\rm Im}\,z\neq 0$, and 
$$\left|\partial_{x'}^\alpha\partial_{\xi'}^\beta\left((c_0r_0-z)^{-1}\right)\right|\le C_{\alpha,\beta}\langle\xi'\rangle^{-2-|\beta|}\eqno{(5.16)}$$
for $|c_0r_0-{\rm Re}\,z|\ge\delta'$. By (5.4), (5.5), (5.6), (5.15) and (5.16),
$$\left|\partial_{x'}^\alpha\partial_{\xi'}^\beta\left(c_1\rho_1-c_2\rho_2)^{-1}\right)\right|\le C_{\alpha,\beta}|{\rm Im}\,z|^{-1-|\alpha|-|\beta|}\eqno{(5.17)}$$
for $(z,x',\xi')\in{\cal M}_1$, $|c_0r_0-{\rm Re}\,z|\le\delta'$,  
$$\left|\partial_{x'}^\alpha\partial_{\xi'}^\beta\left(c_1\rho_1-c_2\rho_2)^{-1}\right)\right|\le C_{\alpha,\beta}|{\rm Im}\,z|^{-\frac{1}{2}-|\alpha|-|\beta|},\quad |\alpha|+|\beta|\ge 1,\eqno{(5.18)}$$
for $(z,x',\xi')\in{\cal M}_1$, $|c_0r_0-{\rm Re}\,z|\ge\delta'$, and 
$$\left|\partial_{x'}^\alpha\partial_{\xi'}^\beta\left((c_1\rho_1-c_2\rho_2)^{-1}\right)\right|\le C_{\alpha,\beta}\langle\xi'\rangle^{-1-|\beta|}\eqno{(5.19)}$$
for $(z,x',\xi')\in{\cal M}_2$. In particular, these estimates imply $|{\rm Im}\,z|(c_1\rho_1-c_2\rho_2)^{-1}\in {\cal S}^{-1}_{\frac{1}{2}-\epsilon}$
for $z\in Z_{1,\epsilon}$ and $ (c_1\rho_1-c_2\rho_2)^{-1}\in {\cal S}^{-1}_0$
for $z\in Z_3$. Hence we have
$$\left\|{\rm Op}_h\left((c_1\rho_1-c_2\rho_2)^{-1}\right)g\right\|_{L^2(\Gamma)}\le \frac{C}{|{\rm Im}\,z|}
\|g\|_{L^2(\Gamma)},\quad\forall g\in L^2(\Gamma).\eqno{(5.20)}$$
By (5.2) and (5.20),
$$\left\|{\rm Op}_h\left((c_1\rho_1-c_2\rho_2)^{-1}\right){\rm Op}_h(c_1\rho_1-c_2\rho_2)f\right\|_{L^2(\Gamma)}\le \frac{Ch}{|{\rm Im}\,z|^{3/2}}
\|f\|_{L^2(\Gamma)}.\eqno{(5.21)}$$
In view of (5.7), (5.8), (5.17)-(5.19), by Proposition 2.2 we have
$$\left\|{\rm Op}_h\left((c_1\rho_1-c_2\rho_2)^{-1}\right){\rm Op}_h(c_1\rho_1-c_2\rho_2)-Id\right\|_{L^2(\Gamma)\to L^2(\Gamma)}\le 
\frac{Ch}{|{\rm Im}\,z|^{5/2}}.\eqno{(5.22)}$$
Combining (5.21) and (5.22) leads to the inequality
$$\left\|f\right\|_{L^2(\Gamma)}\le \frac{Ch}{|{\rm Im}\,z|^{3/2}}
\|f\|_{L^2(\Gamma)}+\frac{Ch}{|{\rm Im}\,z|^{5/2}}\|f\|_{L^2(\Gamma)}.\eqno{(5.23)}$$
 Clearly, it follows from (5.23) that if $h$ is taken small enough, for all $z\in Z_3$ and for 
  $z\in Z_{1,\epsilon}$, $|{\rm Im}\,z|\ge C'h^{2/5}$,  with a sufficiently large constant $C'>0$, we have  $\|f\|_{L^2}=0$, as desired.

Consider now the case $z\in Z_{1,\epsilon}$ under the conditions (1.2), (1.4) and (1.5). It is easy to see that the condition 
(1.5) implies $\frac{c_j}{n_j}|_\Gamma\neq c_0$, $j=1,2$. Hence, if $\delta'>0$ is taken small enough we can arrange that
$|\rho_j|\ge Const>0$ on $|c_0r_0-1|\le \delta'$. Therefore, the functions $a^+=(c_1\rho_1-c_2\rho_2)^{-1}$ and
$a^-=c_1\rho_1-c_2\rho_2$ satisfy (2.4) with $\mu_0=\frac{h}{|{\rm Im}\,z|^2}$, so Proposition 2.2 gives in this case (5.22) (and hence (5.23)) 
with $\frac{h}{|{\rm Im}\,z|^{5/2}}$ replaced by $\frac{h}{|{\rm Im}\,z|^2}$. Thus we obtain that $f\equiv 0$, provided $z\in Z_{1,\epsilon}$
and $h$ taken small enough.

Under the conditions (1.2), (1.4) and (1.6), we have $\rho_1\equiv\rho_2$. Hence (5.18) holds for all $(z,x',\xi')\in{\cal M}_1$,
which again implies (2.4) with $\mu_0=\frac{h}{|{\rm Im}\,z|^2}$, and the desired result follows as above.

It remains to consider the case $z\in Z_2$ under the condition (1.8). It suffices to consider the case $|{\rm Im}\,z|\le \gamma_0$ with some constant $0<\gamma_0\ll 1$, since the case $\gamma_0\le|{\rm Im}\,z|\le 1$ is easy and can be treated as the case $z\in Z_3$ above.
By Proposition 3.6 we have
$$\left\|c_1T_1(z,h)f-c_2T_2(z,h)f\right\|_{L^2(\Gamma)}\le C_Nh^{N-s_d}\left\|f\right\|_{L^2(\Gamma)}\eqno{(5.24)}$$
where $T_j$ is defined by replacing in the definition of the operator $T(z,h)$ from Section 3 the functions $c,n$ by $c_j,n_j$.
Recall that in this case we have $|\rho_j|\ge C\langle\xi'\rangle$, which implies  
$c_1T_1-c_2T_2\in{\rm OP}{\cal S}_0^1$. Since $|c_1\rho_1-c_2\rho_2|\ge C\langle\xi'\rangle$ on supp$\,(1-\chi)$, we have 
$$\left\|{\rm Op}_h(1-\chi)g\right\|_{L^2(\Gamma)}\le C\left\|(c_1T_1-c_2T_2)g\right\|_{L^2(\Gamma)}+O_N(h^N)\left\|g\right\|_{L^2(\Gamma)}\eqno{(5.25)}$$
for every $g\in L^2$ and $N\ge 1$. By (5.24) and (5.25),
$$\left\|{\rm Op}_h(1-\chi)f\right\|_{L^2(\Gamma)}\le C_Nh^{N-s_d}\left\|f\right\|_{L^2(\Gamma)}.\eqno{(5.26)}$$
We will show that
$$|{\rm Im}\,z|\left\|{\rm Op}_h(\chi)f\right\|_{L^2(\Gamma)}^2\le C_Nh^{N-s_d}\left\|f\right\|_{L^2(\Gamma)}^2.\eqno{(5.27)}$$
To this end recall that $z=-1+i{\rm Im}\,z$. Clearly, there exists $0<t\le 1$ so that we can write
$$c_1T_1(z,h)-c_2T_2(z,h)=c_1T_1(-1,h)f-c_2T_2(-1,h)$$ $$+i{\rm Im}\,z\left(c_1\frac{dT_1}{dz}(z_t,h)-c_2\frac{dT_2}{dz}(z_t,h)\right)\eqno{(5.28)}$$
where $z_t=-1+it\,{\rm Im}\,z\in Z_2$. By Lemma 3.9 we have
$$\left|{\rm Re}\,\left\langle (c_1T_1(-1,h)f-c_2T_2(-1,h))f,f\right\rangle_{L^2(\Gamma)}\right|\le C_Nh^{N-s_d}\left\|f\right\|_{L^2(\Gamma)}^2.\eqno{(5.29)}$$
By (5.24), (5.28) and (5.29),
$$|{\rm Im}\,z|\left|{\rm Im}\,\left\langle\left(c_1\frac{dT_1}{dz}(z_t,h)-c_2\frac{dT_2}{dz}(z_t,h)\right)f,f\right\rangle_{L^2(\Gamma)}\right|$$
 $$\le\left|{\rm Re}\,\left\langle (c_1T_1(-1,h)f-c_2T_2(-1,h))f,f\right\rangle_{L^2(\Gamma)}\right|$$ $$+
 \left|{\rm Re}\,\left\langle (c_1T_1(z,h)f-c_2T_2(z,h))f,f\right\rangle_{L^2(\Gamma)}\right|\le C_Nh^{N-s_d}\left\|f\right\|_{L^2(\Gamma)}^2.\eqno{(5.30)}$$
 It follows from (5.30) that to prove (5.27) it suffices to show that
 $$\left\|{\rm Op}_h(\chi)f\right\|_{L^2(\Gamma)}^2\le C\left|{\rm Im}\,\left\langle\left(c_1\frac{dT_1}{dz}(z,h)-c_2\frac{dT_2}{dz}(z,h)\right)f,f\right\rangle_{L^2(\Gamma)}\right|\eqno{(5.31)}$$
 for every $z\in Z_2$ with a constant $C>0$ independent of $z$ and $h$. In view of Lemma 3.9 we have
 $$\left\|\left(c_1\frac{dT_1}{dz}(z,h)-c_2\frac{dT_2}{dz}(z,h)\right)f-{\rm Op}_h(\kappa(z))f\right\|_{L^2(\Gamma)}\le Ch\|f\|_{H^{-1}(\Gamma)}
 \eqno{(5.32)}$$
 where
 $$\kappa(z)=c_1\frac{d\rho_1(z)}{dz}-c_2\frac{d\rho_2(z)}{dz}=-\frac{n_1}{2\rho_1(z)}+\frac{n_2}{2\rho_2(z)}$$
 $$=
 \frac{n_2^2\rho_1^2-n_1^2\rho_2^2}{2\rho_1\rho_2(n_1\rho_2+n_2\rho_1)}= \frac{c_1c_2(n_2^2-n_1^2)r_0-zn_1n_2(c_2n_2-c_1n_1)}{2c_1c_2\rho_1\rho_2(n_1\rho_2+n_2\rho_1)}.$$
 Clearly, we have $\kappa(z)\in{\cal S}^{-1}_0$ and  
 $$\frac{d\kappa(z)}{dz}=-\frac{n_1^2}{4\rho_1(z)^3}+\frac{n_2^2}{4\rho_2(z)^3}=O\left(\langle\xi'\rangle^{-3}\right)$$
 which implies
 $$\kappa(z)=\kappa(-1)+O\left(|{\rm Im}\,z|\langle\xi'\rangle^{-3}\right).\eqno{(5.33)}$$
 Since $\rho_j(-1)=i|\rho_j(-1)|$, we have
 $$i\kappa(-1)=\frac{c_1c_2(n_1^2-n_2^2)r_0+n_1n_2(c_1n_1-c_2n_2)}{2c_1c_2|\rho_1||\rho_2|(n_1|\rho_2|+n_2|\rho_1|)}.$$
 On the other hand, it is easy to see that the condition (1.7) implies
 $$(n_1(x)-n_2(x))(c_1(x)n_1(x)-c_2(x)n_2(x))> 0,\quad\forall x\in\Gamma,$$
 which in turn implies
 $$\left|c_1c_2(n_1^2-n_2^2)r_0+n_1n_2(c_1n_1-c_2n_2)\right|\ge C\langle\xi'\rangle^{2}$$
 and hence
 $$|{\rm Im}\,\kappa(-1)|=|\kappa(-1)|\ge C\langle\xi'\rangle^{-1}.\eqno{(5.34)}$$
 By (5.33) and (5.34),
 $$|{\rm Im}\,\kappa(z)|\ge C\langle\xi'\rangle^{-1}\eqno{(5.35)}$$
 provided $|{\rm Im}\,z|\le \gamma_0$ with some constant $0<\gamma_0\ll 1$. 
 Clearly, we have  
 $${\rm Im}\,\left\langle {\rm Op}_h(\kappa(z))f,f\right\rangle_{L^2(\Gamma)}=\left\langle {\cal A}f,f\right\rangle_{L^2(\Gamma)}$$
 where ${\cal A}=(2i)^{-1}({\rm Op}_h(\kappa(z))-{\rm Op}_h(\kappa(z))^*)$ is an $h$-psdo belonging to ${\rm OP}{\cal S}^{-1}_0$
 with principal symbol ${\rm Im}\,\kappa(z)$. Since the function ${\rm Im}\,\kappa(z)$ is of constant sign, we can use G\"arding's inequality
 together with (5.35) to obtain
 $$\left|\left\langle {\cal A}f,f\right\rangle_{L^2(\Gamma)}\right|\ge C\|f\|^2_{H^{-1/2}(\Gamma)},\quad C>0.\eqno{(5.36)}$$
 Using that $$\left\|{\rm Op}_h(\chi)f\right\|_{L^2(\Gamma)}\le C\|f\|_{H^{-1}(\Gamma)}\le C\|f\|_{H^{-1/2}(\Gamma)}$$
  it is easy to see
 that (5.31) follows from (5.32) and (5.36), provided that $h$ is taken small enough.
 By (5.26) and (5.27) we conclude
 $$|{\rm Im}\,z|\left\|f\right\|_{L^2(\Gamma)}^2\le C_Nh^{N-s_d}\left\|f\right\|_{L^2(\Gamma)}^2.\eqno{(5.37)}$$
If $|{\rm Im}\,z|\ge 2C_Nh^{N-s_d}$, we deduce from (5.37) that $\|f\|_{L^2}=0$. Since $N\gg 1$ is arbitrary, this implies the desired result in this case.
\eproof

G. Vodev

Universit\'e de Nantes,

 D\'epartement de Math\'ematiques, UMR 6629 du CNRS,
 
 2, rue de la Houssini\`ere, BP 92208, 
 
 44332 Nantes Cedex 03, France,
 
 e-mail: vodev@math.univ-nantes.fr


\begin{thebibliography}
\frenchspacing \baselineskip=12 pt plus 1pt minus 1pt



\bibitem{kn:DS} {\sc M. Dimassi and J. Sj\"ostrand}, 
{\em Spectral asymptotics in semi-classical limit}, London Mathematical Society, Lecture Notes Series, {\bf 268}, Cambridge University
Press, 1999.

\bibitem{kn:DP} {\sc M. Dimassi and V. Petkov}, {\em Upper bound for the counting function of interior transmission eigenvalues},
preprint 2013, arXiv: math.SP: 1308.2594v4.

\bibitem{kn:H} {\sc M. Hitrik, K. Krupchyk, P. Ola and L. P\"aiv\"arinta}, {\em The interior transmission problem and bounds
of transmission eigenvalues}, Math. Res. Lett. {\bf 18} (2011), 279-293.

\bibitem{kn:Ho} {\sc L. H\"ormander}, {\em The analysis of linear partial differential operators}, Vol. {\bf 3}, {\em Pseudo-differential
operators}, Springer Verlag, Berlin, 1985. 
 
\bibitem{kn:LV1} {\sc E. Lakshtanov and B. Vainberg}, {\em Remarks on interior transmission 
eigenvalues, Weyl formula and branching billiards}, J. Phys. A: Math. Theor. {\bf 45} (2012), 125202.

\bibitem{kn:LV2} {\sc E. Lakshtanov and B. Vainberg}, {\em Bound on positive interior transmission 
eigenvalues}, Inverse Problems {\bf 28} (2012), 105005.

\bibitem{kn:LV3} {\sc E. Lakshtanov and B. Vainberg}, {\em Application of elliptic theory to the isotropic interior transmission 
eigenvalue problem}, Inverse Problems {\bf 29} (2013), 104003.

\bibitem{kn:LV4} {\sc E. Lakshtanov and B. Vainberg}, {\em Weyl type bound on positive interior transmission 
eigenvalues}, Comm. PDE {\bf  39} (9) (2014),  1729-1740.

\bibitem{kn:PS} {\sc H. Pham and P. Stefanov}, {\em Weyl asymptotics of the transmission eigenvalues for a constant index of
refraction},  Inverse problems and imagining {\bf 8}(3) (2014), 795-810.

\bibitem{kn:PV} {\sc V. Petkov and G. Vodev}, {\em Asymptotics of the number of the interior transmission eigenvalues}, 
J. Spectral Theory, to appear.

\bibitem{kn:R1} {\sc L. Robbiano}, {\em Spectral analysis of interior transmission eigenvalues}, Inverse Problems {\bf 29} (2013), 104001.

\bibitem{kn:R2} {\sc L. Robbiano}, {\em Counting function for interior transmission eigenvalues}, preprint 2013, arXiv: math.AP: 1310.6273.

\bibitem{kn:Sj} {\sc J. Sj\"ostrand}, {\em Singularit\'es analytiques microlocales}, Ast\'erisque, Vol. {\bf 95}, 1982.
{\em }

\bibitem{kn:SjV} {\sc J. Sj\"ostrand and G. Vodev}, {\em Asymptotics of the number of Rayleigh resonances}, Math. Ann. {\bf 309} (1997),
287-306.

\bibitem{kn:SV1} {\sc P. Stefanov and G. Vodev}, {\em Distribution of resonances for the Neumann problem in linear elasticity
outside a strictly convex body}, Duke Math. J. {\bf 78} (1995), 677-714.

\bibitem{kn:SV2} {\sc P. Stefanov and G. Vodev}, {\em Neumann resonances in linear elasticity for an arbitrary body}, Commun. Math. Phys.
{\bf 176} (1996), 645-659.

\bibitem{kn:S} {\sc J. Sylvester}, {\em Transmission eigenvalues in one dimension}, Inverse Problems {\bf 29} (2013), 104009.



\end{thebibliography}
\end{document}